\RequirePackage{fix-cm}
\documentclass[smallextended]{svjour3}       
\smartqed  
\usepackage{graphicx}
\def\Def{\stackrel{\mathrm{def}}{=}}

\def\beq{\begin{equation}}
\def\eeq{\end{equation}}

\def\R{\mathbb{R}}
\def\E{\mathbb{E}}

\def\BI{\begin{itemize}}
\def\EI{\end{itemize}}

\newcommand{\SetEQ}{\setcounter{equation}{0}}
\newcommand{\refLE}[1]{\ensuremath{\stackrel{(\ref{#1})}{\leq}}}
\newcommand{\refEQ}[1]{\ensuremath{\stackrel{(\ref{#1})}{=}}}
\newcommand{\refGE}[1]{\ensuremath{\stackrel{(\ref{#1})}{\geq}}}

\newcommand{\refPE}[1]{\ensuremath{\stackrel{(\ref{#1})}{\preceq}}}
\newcommand{\refSE}[1]{\ensuremath{\stackrel{(\ref{#1})}{\succeq}}}

\usepackage{amssymb}
\usepackage{cite}
\usepackage{float}
\floatstyle{plaintop}
\restylefloat{table}
\usepackage{xcolor}

\usepackage{mathtools}

\usepackage{hyperref}

\numberwithin{equation}{section}
\numberwithin{theorem}{section}
\numberwithin{lemma}{section}
\numberwithin{remark}{section}

\def\Det{{\rm Det}}
\def\Tr{{\rm Tr}}
\def\DFP{{\rm DFP}}
\def\BFGS{{\rm BFGS}}

\def\Broyd{{\rm Broyd}}

\DeclareMathOperator*{\argmax}{argmax}

\renewcommand\arraystretch{1.5}
\def\ba{\begin{array}}
\def\ea{\end{array}}
\def\beann{\begin{eqnarray*}}
\def\eeann{\end{eqnarray*}}
\def\bea{\begin{eqnarray}}
\def\eea{\end{eqnarray}}
\def\bal{\renewcommand\arraystretch{2}\begin{tabular}{|l|}}
\def\eal{\end{tabular}}

\def\BT{\begin{theorem}}
\def\ET{\end{theorem}}
\def\BL{\begin{lemma}}
\def\EL{\end{lemma}}
\def\BC{\begin{corollary}}
\def\EC{\end{corollary}}
\def\BE{\begin{example}}
\def\EE{\end{example}}
\def\BD{\begin{definition}}
\def\ED{\end{definition}}
\def\BR{\begin{remark}}
\def\ER{\end{remark}}
\def\BAS{\begin{assumption}}
\def\EAS{\end{assumption}}
\def\BI{\begin{itemize}}
\def\EI{\end{itemize}}
\def\BCA{\begin{cases}}
\def\ECA{\end{cases}}

\def\BMP{\begin{minipage}{9.5cm}}
\def\EMP{\end{minipage}}
\def\MPT{\begin{minipage}{11.5cm}}
\def\EPT{\end{minipage}}

\def\la{\langle}
\def\ra{\rangle}

\begin{document}

\title{Rates of superlinear convergence for classical
quasi-Newton methods\thanks{Research results in this paper
were obtained with support of ERC Advanced Grant 788368.}}

\author{Anton Rodomanov \and Yurii Nesterov}

\institute{Anton Rodomanov \at
Institute of Information and Communication Technologies,
Electronics and Applied Mathematics, Catholic
University of Louvain, Louvain-la-Neuve, Belgium.\\
\email{anton.rodomanov@uclouvain.be}          
\and
Yurii Nesterov \at
Center for Operations Research and Econometrics,
Catholic University of Louvain, Louvain-la-Neuve,
Belgium.\\
\email{yurii.nesterov@uclouvain.be}
}

\date{%
  Received: 30 March 2020 / Accepted: 18 January 2021 \\
  \textcopyright{} The Author(s) 2021
}

\def\subclassname{{\bfseries Mathematics Subject Classification}\enspace}
\headerboxheight=60pt
\def\makeheadbox{\noindent\small
  Mathematical Programming \\
  \url{https://doi.org/10.1007/s10107-021-01622-5}
}

\maketitle

\begin{abstract}
We study the local convergence of
classical quasi-Newton methods for nonlinear
optimization. Although it was well established a long time
ago that asymptotically these methods converge
superlinearly, the corresponding rates of convergence still
remain unknown. In this paper, we address this problem. We
obtain first explicit non-asymptotic rates of superlinear
convergence for the standard quasi-Newton methods, which are
based on the updating formulas from the convex Broyden
class. In particular, for the well-known DFP and BFGS
methods, we obtain the rates of the form $(\frac{n L^2}
{\mu^2 k})^{k/2}$ and $(\frac{n L}{\mu k})^{k/2}$
respectively, where $k$ is the iteration counter, $n$ is the
dimension of the problem, $\mu$ is the strong convexity
parameter, and $L$ is the Lipschitz constant of the
gradient.
\end{abstract}

\section{Introduction}
\SetEQ

\textbf{Motivation.} In this work, we investigate the
classical quasi-Newton algorithms for smooth unconstrained
optimization, the main examples of which are the
\emph{Davidon--Fletcher--Powell (DFP) method}
\cite{Davidon1959,FletcherPowell1963} and the
\emph{Broyden--Fletcher--Goldfarb--Shanno (BFGS)} method
\cite{Broyden1970p1,Broyden1970p2,Fletcher1970,Goldfarb1970,Shanno1970}.
These algorithms are based on the idea of replacing the
exact Hessian in the Newton method with some approximation,
that is updated in iterations according to certain formulas,
involving only the gradients of the objective function. For
an introduction into the topic, see
\cite{DennisMore1977} and
\cite[Chapter~6]{NocedalWright2006}; also see
\cite{LewisOverton2013} for the treatment of quasi-Newton
algorithms in the context of nonsmooth optimization
and
\cite{GowerGoldfarbRichtarik2016,GowerRichtarik2017,KovalevEtAl2020}
for randomized variants of quasi-Newton methods.

One of the questions about quasi-Newton methods, that has
been extensively studied in the literature, is their
\emph{superlinear convergence}. First theoretical results
here were obtained for the methods with exact line search,
first by Powell \cite{Powell1971}, who analyzed the DFP
method, and then by Dixon \cite{Dixon1972p1,Dixon1972p2},
who showed that with the exact line search all quasi-Newton
algorithms in the Broyden family \cite{Broyden1967}
coincide. Soon after that Broyden, Dennis and Mor{\'e}
\cite{BroydenDennisMore1973} considered the quasi-Newton
algorithms without line search and proved the local
superlinear convergence of DFP, BFGS and several other
methods. Their analysis was based on the Frobenius-norm
potential function. Later, Dennis and
Mor{\'e}~\cite{DennisMore1974} unified the previous proofs
by establishing the necessary and sufficient condition of
superlinear convergence. This condition together with the
original analysis of Broyden, Dennis and Mor{\'e} have been
applied since then in almost every work on quasi-Newton
methods for proving superlinear convergence (see e.g.
\cite{Stachurski1981,GriewankToint1982,EngelsMartinez1991,YabeYamaki1996,WeiYuYuanLian2004,YabeOgasawaraYoshino2007,MokhtariEisenRibeiro2018,GaoGoldfarb2019}).
Finally, one should mention that an important contribution
to the theoretical analysis of quasi-Newton methods has been
made by Byrd, Liu, Nocedal and Yuan in the series of works
\cite{ByrdNocedalYuan1987,ByrdNocedal1989,ByrdLiuNocedal1992},
where they introduced a new potential function by combining
the trace with the logarithm of determinant.

However, the theory of superlinear convergence of
quasi-Newton methods is still far from being complete. The
main reason for this is that all currently existing results
on superlinear convergence of quasi-Newton methods are only
asymptotic: they simply show that the ratio of successive
residuals in the method tends to zero as the number of
iterations goes to infinity, without providing any specific
bounds on the corresponding \emph{rate} of convergence. It
is therefore important to obtain some \emph{explicit} and
\emph{non-asymptotic} rates of superlinear convergence
for quasi-Newton methods.

This observation was the starting point for a recent work
\cite{RodomanovNesterov2020}, where the \emph{greedy}
analogs of the classical quasi-Newton methods have been
developed. As opposed to the classical quasi-Newton methods,
which use the difference of successive iterates for updating
Hessian approximations, these methods employ basis vectors,
greedily selected to maximize a certain measure of progress.
As shown in \cite{RodomanovNesterov2020}, greedy
quasi-Newton methods have superlinear convergence rate
of the form $(1-\frac{\mu}{nL})^{k^2/2} (\frac{n
L}{\mu})^k$, where $k$ is the iteration counter, $n$ is the
dimension of the problem, $\mu$ is the strong convexity
parameter, and $L$ is the Lipschitz constant of the
gradient.

In this work, we continue the same line of research but now
we study the \emph{classical} quasi-Newton methods. Namely,
we consider the methods, based on the updates from the
convex Broyden class, which is formed by all convex
combinations of the DFP and BFGS updates. For this class, we
derive explicit bounds on the rate of superlinear
convergence of standard quasi-Newton methods without line
search. In particular, for the standard DFP and BFGS
methods, we obtain the rates of the form $(\frac{n
L^2}{\mu^2 k})^{k/2}$ and $(\frac{n L}{\mu k})^{k/2}$
respectively.

\textbf{Contents.} This paper is organized as follows.
First, in Section~\ref{sec-broyd-fam}, we study the convex
Broyden class of updating rules for approximating a
self-adjoint positive definite linear operator, and
establish several important properties of this class. Then,
in Section~\ref{sec-quad}, we analyze the standard
quasi-Newton scheme, based on the updating rules from the
convex Broyden class, as applied to minimizing a quadratic
function. We show that this scheme has the same rate of
linear convergence as that of the classical gradient method,
and also a superlinear convergence rate of the form
$(\frac{Q}{k})^{k/2}$, where $Q \geq 1$ is a certain
constant, related to the condition number, and depending on
the method. After that, in Section~\ref{sec-gen}, we
consider the general problem of unconstrained minimization
and the corresponding quasi-Newton scheme for solving it. We
show that, for this scheme, it is possible to prove
absolutely the same results as for the quadratic function,
provided that the starting point is sufficiently close to
the solution. In Section~\ref{sec-disc}, we compare the
rates of superlinear convergence, that we obtain for
the classical quasi-Newton methods, with the corresponding
rates of the greedy quasi-Newton methods.
Section~\ref{sec-app} contains some auxiliary results, that
we use in our analysis.

\textbf{Notation.} In what follows, $\E$ denotes an
arbitrary $n$-dimensional real vector space. Its dual space,
composed by all linear functionals on $\E$, is denoted by
$\E^*$. The value of a linear function $s \in \E^*$,
evaluated at point $x \in \E$, is denoted by $\la s, x
\ra$.

For a smooth function $f : \E \to \R$, we denote by $\nabla
f(x)$ and $\nabla^2 f(x)$ its gradient and Hessian
respectively, evaluated at a point $x \in \E$. Note that
$\nabla f(x) \in \E^*$, and $\nabla^2 f(x)$ is a
self-adjoint linear operator from $\E$ to $\E^*$.

The partial ordering of self-adjoint linear operators is
defined in the standard way. We write $A \preceq A_1$
for $A, A_1 : \E \to \E^*$ if $\la (A_1 - A) x, x \ra
\geq 0$ for all $x \in \E$, and $W \preceq W_1$ for
$W, W_1 : \E^* \to \E$ if $\la s, (W_1 - W) s \ra \geq
0$ for all $s \in \E^*$.

Any self-adjoint positive definite linear operator $A : \E
\to \E^*$ induces in the spaces $\E$ and $\E^*$ the
following pair of conjugate Euclidean norms:
\beq\label{def-norms}
\ba{rclrcl}
\| h \|_A &\Def& \la A h, h \ra^{1/2}, \quad h \in \E,
\qquad\quad
\| s \|_A^* &\Def& \la s, A^{-1} s \ra^{1/2}, \quad s \in
\E^*.
\ea
\eeq
When $A = \nabla^2 f(x)$, where $f : \E \to \R$ is a smooth
function with positive definite Hessian, and $x \in \E$, we
prefer to use notation $\| \cdot \|_x$ and $\| \cdot
\|_x^*$, provided that there is no ambiguity with the
reference function $f$.

Sometimes, in the formulas, involving products of linear
operators, it is convenient to treat $x \in \E$ as a
linear operator from $\R$ to $\E$, defined by $x \alpha =
\alpha x$, and $x^*$ as a linear operator from $\E^*$ to
$\R$, defined by $x^* s = \la s, x \ra$. Likewise, any $s
\in \E^*$ can be treated as a linear operator from $\R$ to
$\E^*$, defined by $s \alpha = \alpha s$, and $s^*$ as a
linear operator from $\E$ to $\R$, defined by $s^* x = \la
s, x \ra$. In this case, $x x^*$ and $s s^*$ are rank-one
self-adjoint linear operators from $\E^*$ to $\E$ and from
$\E^*$ to $\E$ respectively, acting as follows:
$$
\ba{rclrcl}
(x x^*) s &=& \la s, x \ra x,
\qquad
(s s^*) x &=& \la s, x \ra s,
\qquad x \in \E, \ s \in \E^*.
\ea
$$

Given two self-adjoint linear operators $A : \E \to \E^*$
and $W : \E^* \to \E$, we define the trace and the
determinant of $A$ with respect to $W$ as follows:
$$
\ba{rclrcl}
\la W, A \ra &\Def& \Tr(W A), \qquad \Det(W, A) &\Def& \Det
(W A).
\ea
$$
Note that $W A$ is a linear operator from $\E$ to itself,
and hence its trace and determinant are well-defined real
numbers (they coincide with the trace and determinant of the
matrix representation of $W A$ with respect to an arbitrary
chosen basis in the space $\E$, and the result is
independent of the particular choice of the basis). In
particular, if $W$ is positive definite, then $\la W, A \ra$
and $\Det(W, A)$ are respectively the sum and the product of
the eigenvalues\footnote{Recall that, for linear
operators $A, B : \E \to \E^*$, a scalar $\lambda \in \R$ is
called a (relative) eigenvalue of $A$ with respect to $B$ if
$A x = \lambda B x$ for some $x \in \E \setminus \{0\}$. If
$A$, $B$ are self-adjoint and $B$ is positive definite,
it is known that there exist eigenvalues $\lambda_1, \ldots,
\lambda_n \in \R$ and a basis $x_1, \dots, x_n \in \E$, such
that $A x_i = \lambda_i B x_i$, $\| x_i \|_B = 1$, $\la B
x_i, x_j \ra = 0$ for all $1 \leq i, j \leq n$, $i \neq
j$.} of $A$ relative to $W^{-1}$. Observe that $\la \cdot,
\cdot \ra$ is a bilinear form, and for any $x \in \E$, we
have
\beq\label{Axx-tr}
\ba{rcl}
\la A x, x \ra &=& \la x x^*, A \ra.
\ea
\eeq
When $A$ is invertible, we also have
\beq\label{tr-det-A-inv-A}
\ba{rclrcl}
\la A^{-1}, A \ra &=& n, \qquad \Det(A^{-1}, \delta A) &=&
\delta^n.
\ea
\eeq
for any $\delta \in \R$. Also recall the following
multiplicative formula for the determinant:
\beq\label{det-mult}
\ba{rcl}
\Det(W, A) &=& \Det(W, G) \Det(G^{-1}, A),
\ea
\eeq
which is valid for any invertible linear operator $G :
\E \to \E^*$. If the operator $W$ is positive
semidefinite, and $A \preceq A_1$ for some self-adjoint
linear operator $A_1 : \E \to \E^*$, then $\la W, A \ra \leq
\la W, A_1 \ra$ and $\Det(W, A) \leq \Det(W, A_1)$.
Similarly, if $A$ is positive semidefinite and $W \preceq
W_1$ for some self-adjoint linear operator $W_1 : \E^* \to
\E$, then $\la W, A \ra \leq \la W_1, A \ra$ and $\Det(W,
A) \leq \Det(W_1, A)$.

\section{Convex Broyden class}\label{sec-broyd-fam}
\SetEQ

Let $A$ and $G$ be two self-adjoint positive definite linear
operators from $\E$ to $\E^*$, where $A$ is the target
operator, which we want to approximate, and $G$ is the
current approximation of the operator $A$. The \emph{Broyden
family} of quasi-Newton updates of $G$ with respect to $A$
along a direction $u \in \E \setminus \{0\}$, is the
following class of updating formulas, parameterized by a
scalar $\phi \in \R$:
\beq\label{def-broyd}
\ba{rcl}
\Broyd_{\phi}(A, G, u) &\Def& \phi \left[ G - \frac{A u u^*
G + G u u^* A}{\la A u, u \ra} + \left( \frac{\la G u, u
\ra}{\la A u, u \ra} + 1 \right) \frac{A u u^* A}{\la A u, u
\ra} \right] \\
&& + \, (1 - \phi) \left[ G - \frac{G u u^* G}{\la G u, u
\ra} + \frac{A u u^* A}{\la A u, u \ra} \right].
\ea
\eeq
Note that $\Broyd_{\phi}(A, G, u)$ depends on $A$ only
through the product $A u$. For the sake of convenience, we
also define $\Broyd_{\phi}(A, G, u) = G$ when $u = 0$.

Two important members of the Broyden family, DFP and BFGS
updates, correspond to the values $\phi = 1$ and $\phi = 0$
respectively:
\beq\label{def-dfp-bfgs}
\ba{rcl}
\DFP(A, G, u) &\Def& G - \frac{A u u^* G + G u u^* A}{\la
A u, u \ra} + \left( \frac{\la G u, u \ra}{\la A u, u \ra}
+ 1 \right) \frac{A u u^* A}{\la A u, u \ra}, \\
\BFGS(A, G, u) &\Def& G - \frac{G u u^* G}{\la G u, u \ra} +
\frac{A u u^* A}{\la A u, u \ra}.
\ea
\eeq
Thus, the Broyden family consists of all affine
combinations of DFP and BFGS updates:
\beq\label{broyd-aff}
\ba{rcl}
\Broyd_{\phi}(A, G, u) &\refEQ{def-broyd}& \phi \DFP(A, G,
u) + (1 - \phi) \BFGS(A, G, u).
\ea
\eeq
The subclass of the Broyden family, corresponding to $\phi
\in [0, 1]$, is known as the \emph{convex Broyden class} 
(or the \emph{restricted Broyden class} in some texts).

Our subsequent developments will be based on two properties
of the convex Broyden class. The first property states that
each update from this class preserves the bounds on the
relative eigenvalues with respect to the target operator.

\BL\label{lm-pos}
Let $A, G : \E \to \E^*$ be self-adjoint positive definite
linear operators such that
\beq\label{pos-cond}
\ba{rclrcl}
\frac{A}{\xi} &\preceq& G &\preceq& \eta A,
\ea
\eeq
where $\xi, \eta \geq 1$. Then, for any $u \in \E$, and any
$\phi \in [0, 1]$, we have
\beq\label{pos}
\ba{rclrcl}
\frac{A}{\xi} &\preceq& \Broyd_{\phi}(A, G, u) &\preceq&
\eta A.
\ea
\eeq
\EL

\proof
Suppose that $u \neq 0$ since otherwise the claim is
trivial. In view of \eqref{broyd-aff}, it suffices to prove 
\eqref{pos} only for the DFP and BFGS updates independently.

For the DFP update, we have
$$
\ba{rcl}
\DFP(A, G, u) &\refEQ{def-dfp-bfgs}& G - \frac{A u u^* G + G
u u^* A}{\la A u, u \ra} + \left( \frac{\la G u, u \ra}{\la
A u, u \ra} + 1 \right) \frac{A u u^* A}{\la A u, u \ra} \\
&=& \left( I_{\E^*} - \frac{A u u^*}{\la A u, u \ra} \right)
G \left( I_{\E} - \frac{u u^* A}{\la A u, u \ra} \right) +
\frac{A u u^* A}{\la A u, u \ra},
\ea
$$
where $I_{\E}$, $I_{\E^*}$ are the identity operators in the
spaces $\E$, $\E^*$ respectively. Hence,
$$
\ba{rcl}
\DFP(A, G, u) &\refPE{pos-cond}& \eta \left( I_{\E^*} - 
\frac{A u u^*}{\la A u, u \ra} \right) A \left( I_{\E} - 
\frac{u u^* A}{\la A u, u \ra} \right) + \frac{A u u^* A}
{\la A u, u \ra} \\
&=& \eta \left( A - \frac{A u u^* A}{\la A u, u \ra} \right)
+ \frac{A u u^* A}{\la A u, u \ra}
\;=\; \eta A - (\eta - 1) \frac{A u u^* A}{\la A u, u \ra}
\;\preceq\; \eta A, \\
\DFP(A, G, u) &\refSE{pos-cond}& \frac{1}{\xi} \left( I_
{\E^*} - \frac{A u u^*}{\la A u, u \ra} \right) A \left( I_
{\E} - \frac{u u^* A}{\la A u, u \ra} \right) + \frac{A u
u^* A}{\la A u, u \ra} \\
&=& \frac{1}{\xi} \left( A - \frac{A u u^* A}{\la A u, u
\ra} \right) + \frac{A u u^* A}{\la A u, u \ra}
\;=\; \frac{1}{\xi} A + \left( 1 - \frac{1}{\xi} \right) 
\frac{A u u^* A}{\la A u, u \ra}
\;\succeq\; \frac{1}{\xi} A.
\ea
$$

For the BFGS update, we apply Lemma~\ref{lm-proj} (see Appendix):
$$
\ba{rcl}
\BFGS(A, G, u) &\refEQ{def-dfp-bfgs}& G - \frac{G u u^*
G}{\la G u, u \ra} + \frac{A u u^* A}{\la A u, u \ra}
\;\refPE{pos-cond}\; \eta \left( A - \frac{A u u^* A}{\la A
u, u \ra} \right) + \frac{A u u^* A}{\la A u, u \ra} \\ &=&
\eta A - (\eta - 1) \frac{A u u^* A}{\la A u, u \ra}
\;\preceq\; \eta A, \\
\BFGS(A, G, u) &\refEQ{def-dfp-bfgs}& G - \frac{G u u^*
G}{\la G u, u \ra} + \frac{A u u^* A}{\la A u, u \ra}
\;\refSE{pos-cond}\; \frac{1}{\xi} \left( A - \frac{A u u^*
A}{\la A u, u \ra} \right) + \frac{A u u^* A}{\la A u, u
\ra} \\
&=& \frac{1}{\xi} A + \left( 1 - \frac{1}{\xi} \right)
\frac{A u u^* A}{\la A u, u \ra} \;\succeq\; \frac{1}{\xi}
A.$\qed$
\ea
$$

\BR
Lemma~\ref{lm-pos} has first been established in
\cite{Fletcher1970} in a slightly stronger form and using a
different argument. It was also shown there that one of the
relations in \eqref{pos} may no longer be valid if $\phi
\in \R \setminus [0, 1]$.
\ER

The second property of the convex Broyden class, which we
need, is related to the question of convergence of the
approximations $G$ to the target operator~$A$. Note that
without any restrictions on the choice of the update
directions $u$, one cannot guarantee any convergence of $G$
to $A$ in the usual sense (see
\cite{DennisMore1974,RodomanovNesterov2020} for more
details). However, for our goals it will be sufficient to
show that, independently of the choice of $u$, it is still
possible to ensure that $G$ converges to $A$ \emph{along the
update directions $u$}, and estimate the corresponding rate
of convergence.

Let us define the following measure of the
closeness of $G$ to $A$ along the direction $u$:
\beq\label{def-theta}
\ba{rcl}
\theta(A, G, u) &\Def& \left[ \frac{\la (G - A) A^{-1} (G -
A) u, u \ra}{\la G A^{-1} G u, u \ra} \right]^{1/2},
\ea
\eeq
where, for the sake of convenience, we define $\theta(A, G,
u) = 0$ if
$u = 0$. Note that $\theta(A, G, u) = 0$ if and only if $G u
= A u$. Thus, our goal now is to establish some upper bounds
on $\theta$, which will help us to estimate the rate, at
which this measure goes to zero. For this, we will study how
certain potential functions change after one update from the
convex Broyden class, and estimate this change from below by
an appropriate monotonically increasing function of
$\theta$. We will consider two potential functions.

The first one is a simple trace potential function,
that we will use only when we can guarantee that $A
\preceq G$:
\beq\label{def-sigma}
\ba{rclrcl}
\sigma(A, G) &\Def& \la A^{-1}, G - A \ra \;\geq\; 0.
\ea
\eeq

\BL\label{lm-sigma-prog}
Let $A, G : \E \to \E^*$ be self-adjoint positive definite
linear operators such that
\beq\label{sigma-prog-AG}
\ba{rclrcl}
A &\preceq& G &\preceq& \eta A
\ea
\eeq
for some $\eta \geq 1$. Then, for any $\phi \in [0, 1]$ and
any $u \in \E$, we have
\beq\label{sigma-prog}
\ba{rcl}
\sigma(A, G) - \sigma(A, \Broyd_{\phi}(A, G, u)) &\geq&
\left(\phi \frac{1}{\eta} + 1 - \phi \right) \theta^2(A, G,
u).
\ea
\eeq
\EL

\proof
We can assume that $u \neq 0$ since otherwise the claim is
trivial. Denote $G_+ \Def \Broyd_{\phi}(A, G, u)$ and
$\theta \Def \theta(A, G, u)$. Then,
\beq\label{sigma-diff}
\ba{rcl}
&&\sigma(A, G) - \sigma(A, G_+) \;\refEQ{def-sigma}\; \la A^
{-1}, G - G_+ \ra \\
&\refEQ{def-broyd}& 2 \phi \frac{\la G u, u \ra}{\la A u,
u \ra} - \left[ \phi \frac{\la G u, u \ra}{\la A u, u \ra} +
1 \right] + (1 - \phi) \frac{\la G A^{-1} G u, u
\ra}{\la G u, u \ra} \\
&=& \phi \frac{\la G u, u \ra}{\la A u, u \ra}
+ (1 - \phi) \frac{\la G A^{-1} G u, u \ra}{\la G u, u
\ra} - 1 \\
&=& \phi \frac{\la (G - A) u, u \ra}{\la A u, u \ra} + (1 -
\phi) \frac{\la G (A^{-1} - G^{-1}) G u, u \ra}{\la G u, u
\ra}.
\ea
\eeq
Note that
$$
\ba{rclrcl}
0 &\refPE{sigma-prog-AG}& G - A &\refPE{sigma-prog-AG}& 
(\eta - 1) A \;\preceq\; \eta A.
\ea
$$
Therefore\footnote{This is evident when $G - A$ is
non-degenerate. The general case then follows by
continuity.},
\beq\label{GmA-aux}
\ba{rcl}
(G - A) A^{-1} (G - A) &\preceq& \eta (G - A).
\ea
\eeq
Consequently,
\beq\label{sigma-aux-fr1}
\ba{rcl}
\frac{\la (G - A) u, u \ra}{\la A u, u \ra} &\refGE
{GmA-aux}& \frac{1}{\eta} \frac{\la (G - A) A^{-1} (G - A)
u, u \ra}{\la A u, u \ra} \\
&\refGE{sigma-prog-AG}&
\frac{1}{\eta} \frac{\la (G - A) A^{-1} (G - A) u, u \ra}
{\la G A^{-1} G u, u \ra} \;\refEQ{def-theta}\; \frac{1}
{\eta} \theta^2.
\ea
\eeq
At the same time,
$$
\ba{rcl}
(G - A) A^{-1} (G - A) &=& G A^{-1} G - 2 G + A \\
&\refPE{sigma-prog-AG}& G A^{-1} G - G \;=\; G (A^{-1} -
G^{-1}) G.
\ea
$$
Hence,
\beq\label{sigma-aux-fr2}
\ba{rcl}
\frac{\la G (A^{-1} - G^{-1}) G u, u \ra}{\la G u, u \ra}
&\geq& \frac{\la (G - A) A^{-1} (G - A) u, u \ra}{\la G u, u
\ra} \\
&\refGE{sigma-prog-AG}& \frac{\la (G - A) A^{-1} (G - A) u,
u \ra}{\la G A^{-1} G u, u \ra} \;\refEQ{def-theta}\;
\theta^2.
\ea
\eeq
Substituting now \eqref{sigma-aux-fr1} and 
\eqref{sigma-aux-fr2} into \eqref{sigma-diff}, we obtain
\eqref{sigma-prog}.
\qed

The second potential function is more universal since we
can work with it even if the condition $A \preceq G$ is
violated. This function was first introduced in
\cite{ByrdNocedal1989}, and is defined as follows:
\beq\label{def-psi}
\ba{rcl}
\psi(A, G) &\Def& \la A^{-1}, G - A \ra - \ln\Det(A^{-1},
G).
\ea
\eeq
In fact, $\psi$ is nothing else but the Bregman
divergence, generated by the strictly convex function $d(G)
\Def -\ln\Det(B^{-1}, G)$, defined on the set of
self-adjoint positive definite linear operators from $\E$ to
$\E^*$, where $B : \E \to \E^*$ is an arbitrary fixed
self-adjoint positive definite linear operator. Indeed,
$$
\ba{rcl}
\psi(A, G) &\refEQ{tr-det-A-inv-A}& -\ln\Det(B^{-1}, G) +
\ln\Det(B^{-1}, A) - \la -A^{-1}, G - A \ra \\
&=& d(G) - d(A) - \la \nabla d(A), G - A \ra.
\ea
$$
Thus, $\psi(A, G) \geq 0$ and $\psi(A, G) = 0$ if and only
if $G = A$.

Let $\omega : (-1, +\infty) \to \R$ be the univariate
function
\beq\label{def-omega}
\ba{rcl}
\omega(t) &\Def& t - \ln(1 + t) \;\geq\; 0.
\ea
\eeq
Clearly, $\omega$ is a convex function, which is decreasing
on $(-1, 0]$ and increasing on $[0, +\infty)$.
Also, on the latter interval, it satisfies
the following bounds (see \cite[Lemma~5.1.5]
{Nesterov2018Lectures}):
\beq\label{omega-lbd}
\ba{rclrclrcl}
\frac{t^2}{2 (1 + t)} &\leq& \frac{t^2}{2 \left( 1 + 
\frac{2}{3} t \right)} &\leq& \omega(t) &\leq& \frac{t^2}{2
+ t}, \qquad t \geq 0.
\ea
\eeq
Thus, for large values of $t$, the function $\omega
(t)$ is approximately linear in $t$, while for small values
of $t$, it is quadratic.

There is a close relationship between
$\omega$ and the potential function $\psi$. Indeed, if
$\lambda_1, \dots, \lambda_n \geq 0$ are the relative
eigenvalues of $G$ with respect to $A$, then
$$
\ba{rcl}
\psi(A, G) &\refEQ{def-psi}& \sum\limits_{i=1}^n (\lambda_i
- 1 - \ln \lambda_i) \;\refEQ{def-omega}\; \sum\limits_
{i=1}^n \omega(\lambda_i - 1).
\ea
$$

We are going to use the function $\omega$ to estimate from
below the change in the potential function $\psi$, which is
achieved after one update from the convex Broyden class, via
the closeness measure $\theta$. However, first of all, we
need an auxiliary lemma.

\BL\label{lm-aux-omega}
For any real $\alpha \geq \beta > 0$, we have
$$
\ba{rcl}
\alpha - \ln \beta - 1 &\geq& \omega(\sqrt{\alpha \beta -
2 \beta + 1}).
\ea
$$
\EL

\proof
Equivalently, we need to prove that
\beq\label{ineq-omega-2}
\ba{rcl}
\alpha - 1 &\geq& \omega(\sqrt{\alpha \beta - 2 \beta + 1})
+ \ln \beta.
\ea
\eeq

Let us show that the right-hand side of \eqref{ineq-omega-2}
is increasing in $\beta$. This is evident if
$\alpha \geq 2$ because $\omega$ is increasing on $[0,
+\infty)$, so suppose that $\alpha < 2$. Denote
\beq\label{aux-omega-def-delta}
\ba{rcl}
t &\Def& \sqrt{\alpha \beta - 2 \beta + 1} \;=\; \sqrt{1 -
(2 - \alpha) \beta} \;\in\; [0, 1).
\ea
\eeq
Note that $t$ is decreasing in $\beta$. Therefore, it
suffices to prove that the right-hand side of
\eqref{ineq-omega-2} is decreasing in $t$.
But
$$
\ba{rcl}
\omega(\sqrt{\alpha \beta - 2 \beta + 1})
+ \ln \beta &\refEQ{aux-omega-def-delta}& \omega(t) +
\ln \frac{1 - t^2}{2 - \alpha} \\
&=& \omega(t) + \ln(1 - t^2) - \ln(2 - \alpha) \\
&\refEQ{def-omega}& t - \ln(1 + t) + \ln(1 - t^2) - \ln(2 -
\alpha) \\
&=& t + \ln(1 - t) - \ln(2 - \alpha) \\
&\refEQ{def-omega}& -\omega(-t) - \ln(2 - \alpha),
\ea
$$
which is indeed decreasing in $t$ since
$\omega$ is decreasing on $(-1, 0]$.

Thus, it suffices to prove \eqref{ineq-omega-2} only in the
boundary case $\beta = \alpha$:
$$
\ba{rcl}
\alpha - 1 &\geq& \omega(\sqrt{\alpha^2 - 2 \alpha + 1}) +
\ln \alpha \;=\; \omega(|\alpha - 1|) + \ln \alpha,
\ea
$$
or, equivalently, in view of \eqref{def-omega}, that
$$
\ba{rcl}
\omega(\alpha - 1) &\geq& \omega(|\alpha - 1|)
\ea
$$
For $\alpha \geq 1$, this is obvious, so suppose that
$\alpha \leq 1$. It now remains to justify that
\beq\label{omega-neg}
\ba{rcl}
\omega(-t) &\geq& \omega(t),
\ea
\eeq
for all $t \in [0, 1)$. But this easily follows by
integration from the fact that
$$
\ba{rcl}
\frac{d}{dt} \omega(-t) &=& -\omega'(-t) \;\refEQ
{def-omega}\; \frac{t}{1-t} \;\geq\; \frac{t}{1+t} \;\refEQ
{def-omega}\; \omega'(t)
\ea
$$
for all $t \in [0, 1)$.
\qed

Now we are ready to prove the main result.

\BL\label{lm-psi-prog}
Let $A, G : \E \to \E^*$ be self-adjoint positive definite
linear operators such that
\beq\label{psi-prog-AG}
\ba{rclrcl}
\frac{1}{\xi} A &\preceq& G &\preceq& \eta A
\ea
\eeq
for some $\xi, \eta \geq 1$. Then, for any $\phi \in [0, 1]$
and any $u \in \E$, we have
$$
\ba{rcl}
\psi(A, G) - \psi(A, \Broyd_{\phi}(A, G, u)) &\geq& 
\phi \, \omega\left( \frac{\theta(A, G, u)}{\xi^{3/2}
\sqrt{\eta}} \right) + (1-\phi) \omega\left(
\frac{\theta(A, G, u)}{\xi} \right).
\ea
$$
\EL

\proof
Suppose that $u \neq 0$ since otherwise the claim is
trivial. Let us denote $G_+ \Def \Broyd_{\phi}(A, G, u)$ and
$\theta \Def \theta(A, G, u)$. We already know that
$$
\ba{rcl}
\la A^{-1}, G - G_+ \ra &\refEQ{sigma-diff}& \phi 
\frac{\la G u, u \ra}{\la A u, u \ra} + (1-\phi) 
\frac{\la G A^{-1} G u, u \ra}{\la G u, u \ra} - 1.
\ea
$$
Applying now Lemma~\ref{lm-broyd-det-inv}, we obtain
$$
\ba{rcl}
\Det(G^{-1}, G_+) &=& \phi \frac{\la A G^{-1} A u, u
\ra}{\la A u, u \ra} + (1-\phi) \frac{\la A u, u \ra}{\la G
u, u \ra}.
\ea
$$
Thus,
\beq\label{psi-prog-prel}
\ba{rcl}
&&\psi(A, G) - \psi(A, G_+) \\
&\refEQ{def-psi}& \la A^{-1}, G -
G_+ \ra + \ln \Det(A^{-1}, G_+) - \ln \Det(A^{-1}, G) \\
&\refEQ{det-mult}& \la A^{-1}, G - G_+ \ra + \ln \Det(G^
{-1}, G_+) \\
&=& \phi \frac{\la G u, u \ra}{\la A u, u \ra}
+ (1-\phi) \frac{\la G A^{-1} G u, u \ra}{\la G u, u \ra} -
1 + \ln\left[ \phi  \frac{\la A G^{-1} A u, u \ra}{\la A
u, u \ra} + (1-\phi) \frac{\la A u, u \ra} {\la G u, u
\ra} \right] \\
&\geq& \phi \left[ \frac{\la G u, u \ra}{\la
A u, u \ra} + \ln \frac{\la A G^{-1} A u, u \ra}{\la A u, u
\ra} \right] + (1-\phi) \left[ \frac{\la G A^{-1} G u, u
\ra}{\la G u, u \ra} + \ln \frac{\la A u, u \ra}{\la G u, u
\ra} \right] - 1 \\
&=& \phi \left[ \frac{\la G u, u \ra}{\la A u, u \ra}
- \ln \frac{\la A u, u \ra}{\la A G^{-1} A u, u \ra} - 1
\right] + (1-\phi) \left[ \frac{\la G A^{-1} G u, u \ra}{\la
G u, u \ra} - \ln \frac{\la G u, u \ra}{\la A u, u \ra} - 1
\right],
\ea
\eeq
where we have used the concavity of the logarithm.

Denote
\beq\label{def-alpha-beta}
\ba{rcllrcl}
\alpha_1 &\Def& \frac{\la G u, u \ra}{\la A u, u
\ra}, \qquad &\beta_1 &\Def& \frac{\la A u, u \ra}{\la A
G^{-1} A u, u \ra}, \\
\alpha_0 &\Def& \frac{\la G A^{-1} G u, u \ra}{\la G u, u
\ra}, \qquad &\beta_0 &\Def& \frac{\la G u, u \ra}{\la A u,
u \ra}.
\ea
\eeq
Clearly, $\alpha_1 \geq \beta_1$ and $\alpha_0 \geq
\beta_0$ by the Cauchy-Schwartz inequality. Also,
$$
\ba{rcl}
\alpha_1 \beta_1 - 2 \beta_1 + 1 &\refEQ{def-alpha-beta}&
\frac{\la G u, u \ra}{\la A G^{-1} A u, u \ra} - 2 \frac{\la
A u, u \ra}{\la A G^{-1} A u, u \ra} + 1 \;=\; \frac{\la (G
- A) G^{-1} (G - A) u, u \ra}{\la A G^{-1} A u, u \ra} \\
&\refGE{psi-prog-AG}& \frac{1}{\eta} 
\frac{\la (G - A) A^{-1} (G - A) u, u \ra}{\la A G^{-1} A u,
u \ra} \;\refGE{psi-prog-AG}\; \frac{1}{\xi^3 \eta} 
\frac{\la (G - A) A^{-1} (G - A) u, u \ra}{\la G A^{-1} G u,
u \ra} \\
&\refEQ{def-theta}& \frac{\theta^2}{\xi^3 \eta}, \\
\alpha_0 \beta_0 - 2 \beta_0 + 1 &\refEQ{def-alpha-beta}&
\frac{\la G A^{-1} G u, u \ra}{\la A u, u \ra} - 2 \frac{\la
G u, u \ra}{\la A u, u \ra} + 1 \;=\; \frac{\la (G - A) A^
{-1} (G - A) u, u \ra}{\la A u, u \ra} \\
&\refGE{psi-prog-AG}& \frac{1}{\xi^2} \frac{\la (G - A)
A^{-1} (G - A) u, u \ra}{\la G A^{-1} G u, u \ra} \;\refEQ
{def-theta}\; \frac{\theta^2}{\xi^2}.
\ea
$$
Therefore, by Lemma~\ref{lm-aux-omega} and the fact
that $\omega$ is increasing on $[0, +\infty)$, we have
$$
\ba{rcl}
\frac{\la G u, u \ra}{\la A u, u \ra} - \ln
\frac{\la A u, u \ra}{\la A G^{-1} A u, u \ra} - 1 &\geq&
\omega\left( \left[ \frac{\la (G - A) G^{-1} (G - A) u, u
\ra}{\la A G^{-1} A u, u \ra} \right]^{1/2} \right)
\;\geq\; \omega\left( \frac{\theta}{\xi^{3/2} \sqrt{\eta}}
\right), \\
\frac{\la G A^{-1} G u, u \ra}{\la G u, u \ra} - \ln
\frac{\la G u, u \ra}{\la A u, u \ra} - 1 &\geq& \omega\left
( \left[ \frac{\la (G - A) A^{-1} (G - A) u, u \ra}
{\la A u, u \ra} \right]^{1/2} \right) \;\geq\; \omega\left
( \frac{\theta}{\xi} \right).
\ea
$$
Combining these inequalities with \eqref{psi-prog-prel}, we
obtain the claim.
\qed

\section{Unconstrained quadratic minimization}
\label{sec-quad}
\SetEQ

In this section, we study the classical quasi-Newton
methods, based on the updating formulas from the convex
Broyden class, as applied to minimizing the quadratic
function
\beq\label{def-quad}
\ba{rcl}
f(x) &\Def& \frac{1}{2} \la A x, x \ra - \la b, x \ra,
\ea
\eeq
where $A : \E \to \E^*$ is a self-adjoint positive definite
operator, and $b \in \E^*$.

Let $B : \E \to \E^*$ be a self-adjoint positive definite
linear operator, that we will use to initialize our methods.
Denote by $\mu > 0$ the strong convexity parameter of $f$,
and by $L > 0$ the Lipschitz constant of the gradient of
$f$, both measured with respect to $B$:
\beq\label{quad-mu-L}
\ba{rclrcl}
\mu B &\preceq& A &\preceq& L B.
\ea
\eeq

Consider the following standard quasi-Newton scheme for
minimizing \eqref{def-quad}. For the sake of simplicity, we
assume that the constant $L$ is available.
\beq\label{quad-sch-qn}
\bal \hline
\textbf{Initialization:} Choose $x_0 \in \E$. Set $G_0 = L
B$. \\ \hline \textbf{For $k \geq 0$ iterate:} \\
1. Update $x_{k+1} = x_k - G_k^{-1} \nabla f(x_k)$. \\
2. Set $u_k = x_{k+1} - x_k$ and choose $\phi_k \in 
[0, 1]$. \\
3. Compute $G_{k+1} = \Broyd_{\phi_k}(A, G_k, u_k)$. \\
\hline
\eal
\eeq

\BR\label{rm-quad-imp}
In an actual implementation of scheme \eqref{quad-sch-qn},
it is typical to store in memory and update in iterations
the matrix $H_k \Def G_k^{-1}$ instead of $G_k$ (or,
alternatively, the Cholesky decomposition of $G_k$). This
allows one to compute $G_{k+1}^{-1} \nabla f(x_k)$ in
$O(n^2)$ operations. Note that, due to a low-rank structure
of the update \eqref{def-broyd}, $H_k$ can be updated into
$H_{k+1}$ also in $O(n^2)$ operations (for specific
formulas, see e.g. \cite[Section~8]{DennisMore1977}).
\ER

To measure the convergence rate of scheme 
\eqref{quad-sch-qn}, we look at the norm of the gradient,
measured with respect to $A$:
\beq\label{quad-def-lam}
\ba{rcl}
\lambda_f(x) &\Def& \| \nabla f(x) \|_A^* \;\refEQ
{def-norms}\; \la \nabla f(x), A^{-1} \nabla f(x) \ra^{1/2}.
\ea
\eeq

The following lemma shows that the measure $\theta(A, G_k,
u_k)$, that we introduced in \eqref{def-theta} to measure
the closeness of $G_k$ to $A$ along the direction $u_k$, is
directly related to the progress of one step of the scheme
\eqref{quad-sch-qn}. Note that it is important here that the
updating direction $u_k = x_{k+1} - x_k$ is chosen as the
difference of the iterates, and, for other choices of $u_k$,
this result is no longer true.

\BL\label{lm-quad-prog-lam}
In scheme \eqref{quad-sch-qn}, for all $k \geq 0$, we have
\beq\label{quad-lam-via-theta}
\ba{rcl}
\lambda_f(x_{k+1}) &=& \theta(A, G_k, u_k) \lambda_f(x_k).
\ea
\eeq
\EL

\proof
Indeed,
$$
\ba{rcl}
\nabla f(x_{k+1}) &\refEQ{def-quad}& \nabla f(x_k) + A (x_
{k+1} - x_k) \;\refEQ{quad-sch-qn}\; -G_k u_k + A u_k \;=\;
-(G_k - A) u_k.
\ea
$$
Hence, denoting $\theta_k \Def \theta(A, G_k, u_k)$, we
get
$$
\ba{rcl}
\lambda_f(x_{k+1}) &\refEQ{quad-def-lam}& \la (G_k - A)
A^{-1} (G_k - A) u_k, u_k \ra^{1/2}
\;\refEQ{def-theta}\; \theta_k \la G_k A^{-1} G_k u_k, u_k
\ra^{1/2} \\
&\refEQ{quad-sch-qn}& \theta_k \la \nabla f(x_k), A^{-1}
\nabla f(x_k) \ra^{1/2}
\;\refEQ{quad-def-lam}\; \theta_k \lambda_f(x_k).$\qed$
\ea
$$

Let us show that the scheme \eqref{quad-sch-qn} has global
linear convergence, and that the corresponding rate is at
least as good as that of the standard gradient method.

\BT\label{th-quad-lin}
In scheme \eqref{quad-sch-qn}, for all $k \geq 0$, we have
\beq\label{quad-op}
\ba{rclrcl}
A &\preceq& G_k &\preceq& \frac{L}{\mu} A,
\ea
\eeq
and
\beq\label{quad-lin}
\ba{rcl}
\lambda_f(x_k) &\leq& \left( 1 - \frac{\mu}{L} \right)^k
\lambda_f(x_0).
\ea
\eeq
\ET

\proof
For $k=0$, \eqref{quad-op} follows from the fact that $G_0 =
L B$ and \eqref{quad-mu-L}. For all other $k \geq 1$, it
follows by induction using Lemma~\ref{lm-pos}.

Thus, we have
\beq\label{quad-A-m-G}
\ba{rclrcl}
0 &\refPE{quad-op}& A^{-1} - G_k^{-1} &\refPE{quad-op}&
\left( 1 - \frac{\mu}{L} \right) A^{-1}.
\ea
\eeq
Therefore,
$$
\ba{rcl}
(G_k - A) A^{-1} (G_k - A) &=& G_k (A^{-1} - G_k^{-1}) A 
(A^{-1} - G_k^{-1}) G_k \\
&\preceq& \left( 1 - \frac{\mu}{L} \right)^2 G_k A^{-1} G_k,
\ea
$$
and so
$$
\ba{rcl}
\theta(A, G_k, u_k) &\refLE{def-theta}& 1 - \frac{\mu}{L}.
\ea
$$
Applying now Lemma~\ref{quad-lam-via-theta}, we obtain
\eqref{quad-lin}.
\qed

Now, let us establish the superlinear convergence of the
scheme \eqref{quad-sch-qn}. First, we do this by working
with the trace potential function $\sigma$, defined by
\eqref{def-sigma}. Note that this is possible since $A
\preceq G_k$ in view of \eqref{quad-op}.

\BT
In scheme \eqref{quad-sch-qn}, for all $k \geq 1$, we have
\beq\label{quad-super}
\ba{rcl}
\lambda_f(x_k) &\leq& \frac{1}{\prod_{i=0}^{k-1} \left(
\phi_i \frac{\mu}{L} + 1 - \phi_i \right)^{1/2}} \left( 
\frac{n L}{\mu k} \right)^{k/2} \lambda_f(x_0).
\ea
\eeq
\ET

\proof
Denote $\sigma_i \Def \sigma (A, G_i)$, $\theta_i \Def
\theta(A, G_i, u_i)$, and $p_i \Def \phi_i \frac{\mu}{L} + 1
- \phi_i$ for any $i \geq 0$. Let $k \geq 1$ be arbitrary.
From \eqref{quad-op} and Lemma~\ref{lm-sigma-prog}, it
follows that
$$
\ba{rcl}
\sigma_i - \sigma_{i+1} &\geq& p_i \theta_i^2
\ea
$$
for all $0 \leq i \leq k - 1$. Summing up these
inequalities, we obtain
\beq\label{sum-p-theta}
\ba{rcl}
\sum\limits_{i=0}^{k-1} p_i \theta_i^2 &\leq& \sigma_0 -
\sigma_k \;\refLE{def-sigma}\; \sigma_0 \;\refEQ
{quad-sch-qn}\; \sigma(A, L B) \;\refEQ{def-sigma}\; 
\la A^{-1}, L B - A \ra \\
&\refLE{quad-mu-L}& \la A^{-1}, \frac{L}{\mu} A - A \ra
\;\refEQ{tr-det-A-inv-A}\; n \left( \frac{L}{\mu} - 1
\right) \;\leq\; \frac{n L}{\mu}.
\ea
\eeq
Hence, by Lemma~\ref{lm-quad-prog-lam} and the
arithmetic-geometric mean inequality,
$$
\ba{rcl}
\lambda_f(x_k) &=& \lambda_f(x_0) \prod\limits_{i=0}^
{k-1} \theta_i \;=\; \frac{1}{\prod_{i=0}^{k-1} p_i^{1/2}}
\left[ \prod\limits_{i=0}^{k-1} p_i \theta_i^2 \right]^{1/2}
\lambda_f(x_0) \\
&\leq& \frac{1}{\prod_{i=0}^{k-1} p_i^{1/2}} \left( \frac{1}
{k}\sum\limits_{i=0}^{k-1} p_i \theta_i^2 \right)^{k/2}
\lambda_f(x_0) \;\refLE{sum-p-theta}\; \frac{1}{\prod_
{i=0}^{k-1} p_i^{1/2}} \left( \frac{n L}{\mu k} \right)^
{k/2} \lambda_f(x_0).$\qed$
\ea
$$

\BR
As can be seen from \eqref{sum-p-theta}, the factor
$\frac{n L}{\mu}$ in the efficiency estimate 
\eqref{quad-super} can be improved up to $\la A^{-1},
L B - A \ra = \sum_{i=1}^n (\frac{L}{\lambda_i} - 1)$, where
$\lambda_1, \ldots, \lambda_n$ are the eigenvalues of $A$
relative to $B$. This improved factor can be significantly
smaller than the original one if the majority of the
eigenvalues $\lambda_i$ are much larger than $\mu$. However,
for the sake of simplicity, we prefer to work directly with
constants $n$, $L$ and $\mu$. This corresponds to the
worst-case analysis. The same remark applies to all other
theorems on superlinear convergence, that will follow.
\ER

Let us discuss the efficiency estimate \eqref{quad-super}.
Note that its \emph{maximal} value over all $\phi_i \in [0,
1]$ is achieved at $\phi_i = 1$ for all $0 \leq i \leq
k-1$. This corresponds to the \emph{DFP method}. In this
case, the efficiency estimate \eqref{quad-super} looks as
follows:
$$
\ba{rcl}
\lambda_f(x_k) &\leq& \left( \frac{n L^2}{\mu^2 k}
\right)^{k/2} \lambda_f(x_0).
\ea
$$
Hence, the moment, when the superlinear convergence starts,
can be described as follows:
$$
\ba{rclrcl}
\frac{n L^2}{\mu^2 k} &\leq& 1
\qquad \Longleftrightarrow \qquad
k &\geq& \frac{n L^2}{\mu^2}.
\ea
$$
In contrast, the \emph{minimal} value of the efficiency
estimate \eqref{quad-super} over all $\phi_i \in [0, 1]$ is
achieved at $\phi_i = 0$ for all $0 \leq i \leq k-1$. This
corresponds to the \emph{BFGS method}. In this case, the
efficiency estimate \eqref{quad-super} becomes
\beq\label{quad-bfgs-sup}
\ba{rcl}
\lambda_f(x_k) &\leq& \left( \frac{n L}{\mu k} \right)^{k/2}
\lambda_f(x_0),
\ea
\eeq
and the moment, when the superlinear convergence begins,
can be described as follows:
$$
\ba{rclrcl}
\frac{n L}{\mu k} &\leq& 1
\qquad \Longleftrightarrow \qquad
k &\geq& \frac{n L}{\mu}.
\ea
$$
Thus, we see that, compared to DFP, the superlinear
convergence of BFGS starts in $\frac{L}{\mu}$ times earlier,
and its rate is much faster.

Let us present for the scheme \eqref{quad-sch-qn} another
justification of the superlinear convergence rate in the
form \eqref{quad-super}. For this, instead of $\sigma$, we
will work with the potential function~$\omega$, defined by
\eqref{def-omega}. The advantage of this analysis is
that it is extendable onto general nonlinear functions.

\BT
In scheme \eqref{quad-sch-qn}, for all $k \geq 1$, we have
\beq\label{quad-super-2}
\ba{rcl}
\lambda_f(x_k) &\leq& \frac{1}{\prod_{i=0}^{k-1} \left( 
\phi_i \frac{\mu}{L} + 1 - \phi_i \right)^{1/2}} \left( 
\frac{4 n L}{\mu k} \right)^{k/2} \lambda_f(x_0).
\ea
\eeq
\ET

\proof
Denote $\theta_i \Def \theta(A, G_i, u_i)$, $\psi_i \Def
\psi(A, G_i)$, and $p_i \Def \phi_i \frac{\mu}{L} + 1 -
\phi_i$ for any $i \geq 0$. Let $k \geq 1$ and
$0 \leq i \leq k - 1$ be arbitrary. In view of
\eqref{quad-op} and Lemma~\ref{lm-psi-prog}, we have
\beq\label{psi-prog-prel-2}
\ba{rcl}
\psi_i - \psi_{i+1} &\refGE{psi-prog-prel}& \phi_i
\omega\left( \sqrt{\frac{\mu}{L}} \theta_i \right) + 
(1-\phi_i) \omega(\theta_i).
\ea
\eeq
Note that $\theta_i \leq 1$. Indeed, if $u_i = 0$,
then $\theta_i = 0$ by definition. Otherwise,
$$
\ba{rcl}
\theta_i^2 &\refEQ{def-theta}& 1 - \frac{\la (2 G_i - A)
u_i, u_i \ra}{\la G_i A^{-1} G_i u_i, u_i \ra} \;\refLE
{quad-op}\; 1.
\ea
$$
Therefore,
$$
\ba{rcl}
\omega\left( \sqrt{\frac{\mu}{L}} \theta_i \right)
\;\refGE{omega-lbd}\; \frac{\mu}{L} \frac{\theta_i^2}{2
\left( 1 + \sqrt{\frac{\mu}{L}} \theta_i \right)} \;\geq\;
\frac{\mu}{L} \frac{\theta_i^2}{4},
\qquad
\omega(\theta_i) \;\refGE{omega-lbd}\; \frac{\theta_i^2}{2 
(1 + \theta_i)} \;\geq\; \frac{\theta_i^2}{4},
\ea
$$
and we conclude that
$$
\ba{rcl}
\psi_i - \psi_{i+1} &\refGE{psi-prog-prel-2}& \frac{1}{4}
p_i \theta_i^2.
\ea
$$
Summing this inequality and using the fact that $\psi_k \geq
0$, we
obtain
\beq\label{quad-sum-p-theta-omega}
\ba{rcl}
\frac{1}{4} \sum\limits_{i=0}^{k-1} p_i \theta_i^2 &\leq&
\psi_0 - \psi_k \;\leq\; \psi_0 \;\refEQ
{quad-sch-qn}\; \psi(A, L B) \\
&\refEQ{def-psi}& \la A^{-1}, L B - A \ra - \ln \Det(A^{-1},
L B) \\
&\refLE{quad-mu-L}& \la A^{-1}, \frac{L}{\mu} A - A
\ra - \ln \Det(A^{-1}, \frac{L}{\mu} A) \\
&\refEQ{tr-det-A-inv-A}& n \left( \frac{L}{\mu} - 1 - \ln 
\frac{L}{\mu} \right) \;\leq\; \frac{n L}{\mu}.
\ea
\eeq
Hence, by Lemma~\ref{lm-quad-prog-lam} and the
arithmetic-geometric mean inequality,
$$
\ba{rcl}
\lambda_f(x_k) &=& \lambda_f(x_0) \prod\limits_{i=0}^
{k-1} \theta_i \;=\; \frac{1}{\prod_{i=0}^{k-1} p_i^{1/2}}
\left[ \prod\limits_{i=0}^{k-1} p_i \theta_i^2 \right]^{1/2}
\lambda_f(x_0) \\
&\leq& \frac{1}{\prod_{i=0}^{k-1} p_i^{1/2}} \left( \frac{1}
{k} \sum\limits_{i=0}^{k-1} p_i \theta_i^2 \right)^{k/2}
\lambda_f(x_0) \\
&\refLE{quad-sum-p-theta-omega}& \frac{1}
{\prod_{i=0}^{k-1} p_i^{1/2}} \left( \frac{4 n L}{\mu k}
\right)^{k/2} \lambda_f(x_0).$\qed$
\ea
$$

Comparing our new efficiency estimate \eqref{quad-super-2}
with the previous one \eqref{quad-super}, we see that they
differ only in a constant. Thus, for the quadratic
function, we do not gain anything by working with the
potential function $\omega$ instead of $\sigma$.
Nevertheless, our second proof is more universal, and, in
contrast to the first one, can be generalized onto
general nonlinear functions, as we will see in the next
section.

\section{Minimization of general functions}\label{sec-gen}
\SetEQ

Consider now a general unconstrained minimization problem:
\beq\label{prob-gen}
\ba{rcl}
\min\limits_{x \in \E} f(x),
\ea
\eeq
where $f : \E \to \R$ is a twice differentiable function
with positive definite Hessian.

To write down the standard quasi-Newton scheme for
\eqref{prob-gen}, we fix some self-adjoint positive definite
linear operator $B : \E \to \E^*$ and a constant $L > 0$,
that we use to define the initial Hessian approximation.
\beq\label{sch-qn}
\bal \hline
\textbf{Initialization:} Choose $x_0 \in \E$. Set $G_0 = L
B$. \\ \hline
\textbf{For $k \geq 0$ iterate:} \\
1. Update $x_{k+1} = x_k - G_k^{-1} \nabla f(x_k)$. \\
2. Set $u_k = x_{k+1} - x_k$ and choose $\phi_k \in [0, 1]$.
\\
3. Denote $J_k = \int_0^1 \nabla^2 f(x_k + t u_k) dt$.\\
4. Set $G_{k+1} = \Broyd_{\phi_k}(J_k, G_k, u_k)$. \\
\hline
\eal
\eeq

\BR
Similarly to Remark~\ref{rm-quad-imp}, when
implementing scheme~\eqref{sch-qn}, it is common to work
directly with the inverse $H_k \Def G_k^{-1}$ instead of
$G_k$. Also note that it is not necessary to compute $J_k$
explicitly. Indeed, for implementing the Hessian
approximation update at Step~4 (or the corresponding update
for its inverse), one only needs the product
$$
\ba{rcl}
J_k u_k &=& \nabla f(x_{k+1}) - \nabla f(x_k),
\ea
$$
which is just the difference of the successive gradients.
\ER

In what follows, we make the following assumptions about the
problem~\eqref{prob-gen}. First, we assume that, with
respect to the operator $B$, the objective function $f$ is
\emph{strongly convex} with parameter $\mu > 0$ and its
gradient is \emph{Lipschitz continuous} with constant $L$,
i.e.
\beq\label{mu-L}
\ba{rclrcl}
\mu B &\preceq& \nabla^2 f(x) &\preceq& L B
\ea
\eeq
for all $x \in \E$. Second, we assume that the objective
function $f$ is \emph{strongly self-concordant} with some
constant $M \geq 0$, i.e.
\beq\label{sscf}
\ba{rcl}
\nabla^2 f(y) - \nabla^2 f(x) &\preceq& M \| y - x \|_z
\nabla^2 f(w)
\ea
\eeq
for all $x, y, z, w \in \E$. The class of strongly
self-concordant functions was recently introduced in 
\cite{RodomanovNesterov2020}, and contains at least all
strongly convex functions with Lipschitz continuous Hessian 
(see \cite[Example~4.1]{RodomanovNesterov2020}).
It gives us the the following convenient relations between
the Hessians of the objective function:
\BL[see {\cite[Lemma 4.1]{RodomanovNesterov2020}}]
\label{lm-hess-scf}
Let $x, y \in \E$, and let $r \Def \| y - x \|_x$. Then,
\beq\label{hess-x-y}
\ba{rclrcl}
\frac{\nabla^2 f(x)}{1 + M r} &\preceq& \nabla^2 f(y)
&\preceq& (1 + M r) \nabla^2 f(x).
\ea
\eeq
Also, for $J \Def \int_0^1 \nabla^2 f(x + t (y - x)) dt$, we
have
\beq\label{hess-J-x}
\ba{rclrcl}
\frac{\nabla^2 f(x)}{1 + \frac{M r}{2}} &\preceq& J
&\preceq& \left( 1 + \frac{M r}{2} \right) \nabla^2 f(x),
\ea
\eeq
\beq\label{hess-J-y}
\ba{rclrcl}
\frac{\nabla^2 f(y)}{1 + \frac{M r}{2}} &\preceq& J
&\preceq& \left( 1 + \frac{M r}{2} \right) \nabla^2 f(y).
\ea
\eeq
\EL

As a particular example of a nonquadratic function,
satisfying assumptions~\eqref{mu-L}, \eqref{sscf}, one can
consider the regularized log-sum-exp function, defined by
$f(x) \Def \ln(\sum_{i=1}^m e^{\la a_i, x \ra + b_i}) +
\frac{\mu}{2} \| x \|^2$, where $a_i \in \E^*$, $b_i \in \R$
for $i = 1, \ldots, m$, and $\mu > 0$, $\| x \| \Def \la B
x, x \ra^{1/2}$.

\BR
Since we are interested in \emph{local} convergence, it is
possible to relax our assumptions by requiring that
\eqref{mu-L}, \eqref{sscf} hold only in some neighborhood of
a minimizer $x^*$. For this, it suffices
to assume that the Hessian of $f$ is Lipschitz continuous in
this neighborhood with $\nabla^2 f(x^*)$ being positive
definite. These are exactly the standard assumptions, used
in \cite{DennisMore1977} and many other works, studying
local convergence of quasi-Newton methods. However, to
avoid excessive technicalities, we do not do this.
\ER

Let us now analyze the process \eqref{sch-qn}. For measuring
its convergence, we look at the local norm of the gradient:
\beq\label{def-lam}
\ba{rcl}
\lambda_f(x) &\Def& \| \nabla f(x) \|_x^* \;\refEQ
{def-norms}\; \la \nabla f(x), \nabla^2 f(x)^{-1} \nabla f
(x) \ra^{1/2}, \qquad x \in \E.
\ea
\eeq

First, let us estimate the progress of one step of the
scheme \eqref{sch-qn}. Recall that $\theta(J_k, G_k, u_k)$
is the measure of closeness of $G_k$ to $J_k$ along the
direction $u_k$ (see \eqref{def-theta}).

\BL\label{lm-lam-prog}
In scheme \eqref{sch-qn}, for all $k \geq 0$ and $r_k
\Def \| u_k \|_{x_k}$, we have
$$
\ba{rcl}
\lambda_f(x_{k+1}) &\leq& \left( 1 + \frac{M r_k}{2} \right)
\theta(J_k, G_k, u_k) \lambda_f(x_k).
\ea
$$
\EL

\proof
Denote $\theta_k \Def \theta(J_k, G_k, u_k)$. In view of
Taylor's formula,
\beq\label{next-grad}
\ba{rcl}
\nabla f(x_{k+1}) &=& \nabla f(x_k) + J_k (x_{k+1} - x_k)
\;\refEQ{sch-qn}\; -(G_k - J_k) u_k.
\ea
\eeq
Therefore,
$$
\ba{rcl}
\lambda_f(x_{k+1}) &\refEQ{def-lam}& \la \nabla f(x_{k+1}),
\nabla^2 f(x_{k+1})^{-1} \nabla f(x_{k+1}) \ra^{1/2} \\
&\refLE{hess-J-y}& \sqrt{1 + \frac{M r_k}{2}} \, \la \nabla f
(x_{k+1}), J_k^{-1} \nabla f (x_{k+1}) \ra^{1/2} \\
&\refEQ{next-grad}& \sqrt{1 + \frac{M r_k}{2}} \, \la (G_k
- J_k) J_k^{-1} (G_k - J_k) u_k, u_k \ra^{1/2} \\
&\refEQ{def-theta}& \sqrt{1 + \frac{M r_k}{2}} \, \theta_k
\la G_k J_k^{-1} G_k u_k, u_k \ra^{1/2} \\
&\refEQ{sch-qn}& \sqrt{1 + \frac{M r_k}{2}} \, \theta_k \la
\nabla f(x_k), J_k^{-1} \nabla f(x_k) \ra^{1/2} \\
&\refLE{hess-J-x}& \left( 1 + \frac{M r_k}{2} \right) \,
\theta_k \la \nabla f(x_k), \nabla^2 f(x_k)^{-1} \nabla f
(x_k) \ra^{1/2} \\ &\refEQ{def-lam}& \left( 1 + \frac{M
r_k}{2} \right) \theta_k \lambda_f(x_k).$\qed$
\ea
$$

Our next result states that, if the
starting point in scheme \eqref{sch-qn} is chosen
sufficiently close to the solution, then the relative
eigenvalues of the Hessian approximations $G_k$ with respect
to both the Hessians $\nabla^2 f(x_k)$ and the integral
Hessians $J_k$ are always located between 1 and
$\frac{L}{\mu}$, up to some small numerical constant. As a
consequence, the process~\eqref{sch-qn} has at least the
linear convergence rate of the gradient method.

\BT
Suppose that, in scheme \eqref{sch-qn},
\beq\label{lam-ini}
\ba{rcl}
M \lambda_f(x_0) &\leq& \frac{\ln \frac{3}{2}}{4} \frac{\mu}
{L}.
\ea
\eeq
Then, for all $k \geq 0$, we have
\beq\label{op-lin-hess}
\ba{rclrcl}
\frac{1}{\xi_k} \nabla^2 f(x_k) &\preceq& G_k &\preceq&
\xi_k \frac{L}{\mu} \nabla^2 f(x_k),
\ea
\eeq
\beq\label{op-lin-J}
\ba{rclrcl}
\frac{1}{\xi_k'} J_k &\preceq& G_k &\preceq& \xi_k' \frac{L}
{\mu} J_k,
\ea
\eeq
\beq\label{lam-lin}
\ba{rcl}
\xi_k \lambda_f(x_k) &\leq& \left(1 - \frac{\mu}{2 L}
\right)^k \lambda_f(x_0),
\ea
\eeq
where\footnote{Here we follow the standard convention that
the sum over the empty set is defined as zero. Thus, $\xi_0
= e^0 = 1$.}
\beq\label{def-xi}
\ba{rclrclrclrcl}
\xi_k &\Def& e^{M \sum_{i=0}^{k-1} r_i} &\leq& \left( 1 +
\frac{M r_k}{2} \right) e^{M \sum_{i=0}^{k-1} r_i} &\Def&
\xi_k'&\leq& \sqrt{\frac{3}{2}},
\ea
\eeq
and $r_i \Def \| u_i \|_{x_i}$ for any $i \geq 0$.
\ET

\proof
Note that $\xi_0 = 1$ and $G_0 = L B$. Therefore, for $k =
0$, both \eqref{op-lin-hess}, \eqref{lam-lin} are
satisfied. Indeed, the first one reads $\nabla^2 f(x_0)
\preceq L B \preceq \frac{L}{\mu} \nabla^2 f(x_0)$ and
follows from~\eqref{mu-L}, while the second one reads
$\lambda_f(x_0) \leq \lambda_f(x_0)$ and is obviously true.

Now assume that $k \geq 0$, and that \eqref{op-lin-hess},
\eqref{lam-lin} have already been proved for all $0 \leq k'
\leq k$. Combining \eqref{op-lin-hess} with 
\eqref{hess-J-x}, using the definition of $\xi_k'$, we
obtain \eqref{op-lin-J}. Further, denote $\lambda_i \Def
\lambda_f(x_i)$ for $0 \leq i \leq k$. Note that
\beq\label{r-via-lam}
\ba{rcl}
r_k &\refEQ{sch-qn}& \| G_k^{-1} \nabla f(x_k) \|_{x_k}
\;\refEQ{def-norms}\; \la \nabla f(x_k), G_k^{-1} \nabla^2 f
(x_k) G_k^{-1} \nabla f(x_k) \ra^{1/2} \\ &\refLE
{op-lin-hess}& \xi_k \la \nabla f(x_k), \nabla^2 f(x_k)^{-1}
\nabla f(x_k) \ra^{1/2} \;\refEQ{def-lam}\; \xi_k \lambda_k.
\ea
\eeq
Therefore,
\beq\label{sum-r}
\ba{rcl}
M \sum\limits_{i=0}^k r_i &\refLE{r-via-lam}& M \sum\limits_
{i=0}^k \xi_i \lambda_i \;\refLE{lam-lin}\; M \lambda_0
\sum\limits_{i=0}^k \left( 1 - \frac{\mu}{2 L} \right)^i \\
&\leq& \frac{2 L}{\mu} M \lambda_0 \;\refLE{lam-ini}\; 
\frac{\ln \frac{3}{2}}{2}.
\ea
\eeq
Consequently, by the definition of $\xi_k$ and $\xi_k'$,
$$
\ba{rcl}
\xi_k &\leq& \xi_k' \;\leq\; e^{\frac{M r_k}{2}} e^{M \sum_
{i=0}^{k-1} r_i} \;\leq\; e^{M \sum_{i=0}^k r_i}
\;\refLE{sum-r}\; \sqrt{\frac{3}{2}}.
\ea
$$
Thus, \eqref{op-lin-J}, \eqref{def-xi} are now proved. To
finish the proof by induction, it remains to prove
\eqref{op-lin-hess}, \eqref{lam-lin} for $k' = k + 1$.

We start with \eqref{op-lin-hess}. Applying
Lemma~\ref{lm-pos}, using \eqref{op-lin-J}, we obtain
\beq\label{Gp-J}
\ba{rclrcl}
\frac{1}{\xi_k'} J_k &\preceq& G_{k+1} &\preceq& \xi_k'
\frac{L}{\mu} J_k.
\ea
\eeq
Consequently,
$$
\ba{rcl}
G_{k+1} &\refPE{hess-J-y}& \left( 1 + \frac{M r_k}{2}
\right) \xi_k' \frac{L}{\mu} \nabla^2 f(x_{k+1}) \;\refEQ
{def-xi}\; \left(1 + \frac{M r_k}{2} \right)^2 \xi_k 
\frac{L}{\mu} \nabla f(x_{k+1}) \\
&\preceq& e^{M r_k} \xi_k \frac{L}{\mu} \nabla^2 f(x_{k+1})
\;\refEQ{def-xi}\; \xi_{k+1} \frac{L}{\mu} \nabla^2 f(x_
{k+1}),
\ea
$$
and
$$
\ba{rcl}
G_{k+1} &\refSE{hess-J-y}& \frac{\nabla^2 f(x_{k+1})}{\left(
1 + \frac{M r_k}{2} \right) \xi_k'} \;\refEQ{def-xi}\;
\frac{\nabla^2 f (x_{k+1})}{\left( 1 + \frac{M r_k}{2}
\right)^2 \xi_k} \;\succeq\;
\frac{\nabla^2 f(x_{k+1})}{e^{M r_k} \cdot \xi_k}
\;\refEQ{def-xi}\; \frac{\nabla^2 f(x_{k+1})}{\xi_
{k+1}}.
\ea
$$
Thus, \eqref{op-lin-hess} is proved for $k'= k+1$.

It remains to prove \eqref{lam-lin} for $k'= k+1$. By
Lemma~\ref{lm-lam-prog},
\beq\label{lam-next-prel}
\ba{rcl}
\lambda_{k+1} &\leq& \left( 1 + \frac{M r_k}{2} \right)
\theta_k \lambda_k,
\ea
\eeq
where $\theta_k \Def \theta(J_k, G_k, u_k)$. Note that
$$
\ba{rclrcl}
-\left( 1 - \frac{\mu}{\xi_k' L} \right) J_k^{-1}
&\refPE{op-lin-J}& G_k^{-1} - J_k^{-1} &\refPE{op-lin-J}&
(\xi_k' - 1) J_k^{-1}.
\ea
$$
Hence,
$$
\ba{rcl}
(J_k^{-1} - G_k^{-1}) J_k (J_k^{-1} - G_k^{-1}) &\preceq&
\rho_k^2 J_k^{-1},
\ea
$$
where
\beq\label{def-rho}
\ba{rcl}
\rho_k &\Def& \max\left\{ 1 - \frac{\mu}{\xi_k' L} , \,
\xi_k' - 1 \right\} \;\refGE{def-xi}\; 0.
\ea
\eeq
Therefore,
$$
\ba{rcl}
\theta_k^2 &\refEQ{def-theta}& \frac{\la (J_k - G_k)
J_k^{-1} (J_k - G_k) u_k, u_k \ra}{\la J_k G_k^{-1} J_k u_k,
u_k \ra} \;=\; \frac{\la G_k (J_k^{-1} - G_k^{-1}) J_k (J_k^
{-1} - G_k^{-1}) G_k u_k, u_k \ra}{\la J_k G_k^{-1} J_k u_k,
u_k \ra} \;\leq\; \rho_k^2.
\ea
$$
Thus,
$$
\ba{rcl}
\lambda_{k+1} &\refLE{lam-next-prel}& \left( 1 + \frac{M
r_k}{2} \right) \rho_k \lambda_k.
\ea
$$
Consequently,
$$
\ba{rcl}
\xi_{k+1} \lambda_{k+1} &\leq& \xi_{k+1} \left( 1 + \frac{M
r_k}{2} \right) \rho_k \lambda_k \;\refEQ{def-xi}\; e^
{M r_k} \left( 1 + \frac{M r_k}{2} \right) \rho_k \xi_k
\lambda_k \\
&\leq& e^{\frac{3 M r_k}{2}} \rho_k \xi_k \lambda_k \;\refLE
{lam-lin}\; e^{\frac{3 M r_k}{2}} \rho_k \left( 1 - 
\frac{\mu}{2 L} \right)^k \lambda_0.
\ea
$$
It remains to show that
\beq\label{e-rho}
\ba{rcl}
e^{\frac{3 M r_k}{2}} \rho_k &\leq& 1 - \frac{\mu}{2 L}.
\ea
\eeq
Note that
\beq\label{def-zeta}
\ba{rcl}
\zeta_k &\Def& \frac{3 M r_k}{2} \;\refLE{r-via-lam}\; 
\frac{3 M \xi_k \lambda_k}{2} \;\refLE{lam-lin}\; \frac{3 M
\lambda_0}{2} \\
&\refLE{lam-ini}& \frac{3 \ln \frac{3}{2}}{8} \frac{\mu}{L}
\;\leq\; \frac{3 \mu}{16 L} \;\leq\; \frac{\mu}{5 L} \;\leq\;
\frac{1}{5}.
\ea
\eeq
Hence,
\beq\label{lin-aux1}
\ba{rcl}
e^{\zeta_k} &\leq& \sum\limits_{i=0}^{\infty} \zeta_k^i =
1 + \zeta_k \sum\limits_{i=0}^\infty \zeta_k^i 
\;\refLE{def-zeta}\; 1 + \zeta_k \sum\limits_{i=0}^\infty
\left( \frac{1}{5}\right)^i \\
&=& 1 + \frac{5 \zeta_k}{4} \;\refLE{def-zeta}\; 1 + 
\frac{\mu}{4 L}.
\ea
\eeq
Also,
\beq\label{lin-aux2}
\ba{rcl}
\xi_k' &\refLE{def-xi}& \sqrt{\frac{3}{2}} \;\leq\;
\frac{4}{3}.
\ea
\eeq
Combining \eqref{lin-aux1} and \eqref{lin-aux2}, we obtain 
$$
\ba{rcl}
e^{\frac{3 M r_k}{2}} \left( 1 - \frac{\mu}{\xi_k' L}
\right) &\leq& \left( 1 + \frac{\mu}{4 L} \right)
\left( 1 - \frac{3 \mu}{4 L} \right) \;\leq\; 1 - \left( 
\frac{3}{4} - \frac{1}{4} \right) \frac{\mu}{L} \;=\; 1 -
\frac{\mu}{2 L},
\ea
$$
and
$$
\ba{rcl}
e^{\frac{3 M r_k}{2}} (\xi_k' - 1) &\leq& \left( 1 + 
\frac{1}{4} \right) \left( \sqrt{\frac{3}{2}} - 1 \right)
\;=\; \frac{\frac{5}{4} \cdot \frac{1}{2}}{\sqrt{\frac{3}{2}}
+ 1} \;\leq\; \frac{5}{16} \;\leq\; \frac{1}{2} \;\leq\; 1 -
\frac{\mu}{2 L}.
\ea
$$
Thus,
$$
\ba{rcl}
e^{\frac{3 M r_k}{2}} \rho_k &\refEQ{def-rho}& e^{\frac{3 M
r_k}{2}} \max\left\{ 1 - \frac{\mu}{\xi_k' L}, \,
\xi_k' - 1 \right\} \;\leq\; 1 - \frac{\mu}{2 L},
\ea
$$
and \eqref{e-rho} follows.
\qed

Now we are ready to prove the main result of this section
on the superlinear convergence of the scheme
\eqref{sch-qn}. In contrast to the quadratic case, now we
cannot use the proof, based on the trace potential function
$\sigma$, defined by~\eqref{def-sigma}, because we cannot
longer guarantee that $J_k \preceq G_k$. However, the proof,
based on the potential function $\psi$, defined by
\eqref{def-psi}, still works.

\BT
Suppose that the initial point $x_0$ in
scheme~\eqref{sch-qn} is chosen sufficiently close to the
solution, as specified by \eqref{lam-ini}. Then, for all $k
\geq 1$, we have
$$
\ba{rcl}
\lambda_f(x_k) &\leq& \frac{1}{\prod_{i=0}^{k-1} \left(
\phi_i \frac{\mu}{L} + 1 - \phi_i \right)^{1/2}} \left( 
\frac{11 n L} {\mu k} \right)^{k/2} \lambda_f(x_0).
\ea
$$
\ET

\proof
Denote $r_i \Def \| u_i \|_{x_i}$, $\theta_i \Def
\theta(J_i, G_i, u_i)$, $\psi_i \Def \psi(J_i, G_i)$,
$\tilde{\psi}_{i+1} \Def \psi(J_i,G_{i+1})$, and $p_i \Def
\phi_i \frac{\mu}{L} + 1 - \phi_i$ for any $i \geq 0$. Let
$k \geq 1$ and $0 \leq i \leq k - 1$ be arbitrary. By 
\eqref{op-lin-J}, \eqref{def-xi} and
Lemma~\ref{lm-psi-prog}, we have
\beq\label{gen-psi-diff-prel}
\ba{rcl}
\psi_i - \tilde{\psi}_{i+1} &\geq& \phi_i \omega\left(
\frac{2}{3} \sqrt{\frac{\mu}{L}} \theta_i \right) +
(1-\phi_i) \omega\left( \sqrt{\frac{2}{3}} \theta_i \right).
\ea
\eeq
Moreover, since
$$
\ba{rcl}
\theta_i^2 &\refEQ{def-theta}& \frac{\la (G_i - J_i) J_i^
{-1} (G_i - J_i) u_i, u_i \ra}{\la G_i J_i^{-1} G_i u_i,
u_i \ra} \;=\; 1 - \frac{\la (2 G_i - J_i) u_i, u_i \ra}
{\la G_i J_i^{-1} G_i u_i, u_i \ra} \;\refLE{op-lin-J}\; 1,
\ea
$$
we also have
$$
\ba{rcl}
\omega\left( \frac{2}{3} \sqrt{\frac{\mu}{L}} \theta_i
\right) &\refGE{omega-lbd}& \frac{\frac{4}{9} \frac{\mu}
{L} \theta_i^2}{2 \left( 1 + \frac{2}{3} \cdot \frac{2}{3} 
\sqrt{\frac{\mu}{L}} \theta_i \right)} \;\geq\; \frac{
\frac{4}{9}}{2 (1 + \frac{4}{9})} \frac{\mu}{L} \theta_i^2
\;=\; \frac{2}{13} \frac{\mu}{L} \theta_i^2 \;\geq\; 
\frac{1}{7} \frac{\mu}{L} \theta_i^2, \\
\omega\left( \sqrt{\frac{2}{3}} \theta_i \right) &\refGE
{omega-lbd}& \frac{\frac{2}{3} \theta_i^2}{2 \left( 1 +
\sqrt{\frac{2}{3}} \theta_i \right)} \;\geq\; \frac{
\frac{2}{3}}{2 \left( 1 + \sqrt{\frac{2}{3}} \right)}
\;\geq\; \frac{\frac{2}{3}}{4} \theta_i^2 \;=\; 
\frac{1}{6} \theta_i^2 \;\geq\; \frac{1}{7} \theta_i^2.
\ea
$$

Thus,
\beq\label{p-theta-prel}
\ba{rcl}
\frac{1}{7} p_i \theta_i^2 &\refLE{gen-psi-diff-prel}&
\psi_i - \tilde{\psi}_{i+1} \;=\; \psi_i - \psi_{i+1} +
\Delta_i,
\ea
\eeq
where
\beq\label{def-Delta}
\ba{rcl}
\Delta_i \Def \psi_{i+1} - \tilde{\psi}_{i+1} &\refEQ
{def-psi}& \la J_{i+1}^{-1} - J_i^{-1}, G_{i+1} \ra + \ln
\Det (J_i^{-1}, J_{i+1}).
\ea
\eeq
Let us estimate $\sum_{i=0}^{k-1} \Delta_i$ from above. Note
that
\beq\label{J-Jp}
\ba{rcl}
J_{i+1} &\refSE{hess-J-x}& \frac{\nabla^2 f(x_{i+1})}{1 +
\frac{M r_{i+1}}{2}} \;\refSE{hess-J-y}\; \frac{1}{\delta_i}
J_i,
\ea
\eeq
where
\beq\label{def-delta}
\ba{rcl}
\delta_i &\Def& \left( 1 + \frac{M r_{i+1}}{2} \right)
\left( 1 + \frac{M r_i}{2} \right).
\ea
\eeq
Hence,
$$
\ba{rcl}
\la J_{i+1}^{-1} - J_i^{-1}, G_{i+1} \ra &\refLE{J-Jp}& (1 -
\delta_i^{-1}) \la J_{i+1}^{-1}, G_{i+1} \ra \\
&\refLE{op-lin-J}& (1 - \delta_i^{-1}) \sqrt{\frac{3}{2}} 
\frac{L}{\mu} \la J_{i+1}^{-1}, J_{i+1} \ra \\
&\refEQ{tr-det-A-inv-A}& \sqrt{\frac{3}{2}} \frac{n L}{\mu}
(1 - \delta_i^{-1}) \;\refLE {lin-aux2}\; \frac{4 n L}{3
\mu} (1 - \delta_i^{-1}),
\ea
$$
and
\beq\label{sum-Delta-prel}
\ba{rcl}
\sum\limits_{i=0}^{k-1} \Delta_i
&\refLE{def-Delta}& \frac{4 n L}{3 \mu} \sum\limits_{i=0}^
{k-1} (1 - \delta_i^{-1}) + \sum\limits_{i=0}^{k-1} \ln \Det
(J_i^{-1}, J_{i+1}) \\ &=& \frac{4 n L}{3
\mu}\sum\limits_{i=0}^ {k-1}(1 - \delta_i^{-1}) + \ln\Det
(J_0^{-1}, J_k).
\ea
\eeq
At the same time,
$$
\ba{rcl}
\sum\limits_{i=0}^{k-1} (1 - \delta_i^{-1}) &\leq&
\sum\limits_{i=0}^{k-1} \left( 1 - e^{-\frac{M (r_i +
r_{i+1})}{2}} \right) \;\leq\; \frac{M}{2}
\sum\limits_{i=0}^{k-1} (r_i + r_{i+1}) \;\leq\; M
\sum\limits_{i=0}^k r_i \\
&\refLE{lam-lin}& M \lambda_0 \sum\limits_{i=0}^{k-1} \left(
1 - \frac{\mu}{2 L} \right)^i \;\leq\; \frac{2 L}{\mu} M
\lambda_0 \;\refLE{lam-ini}\; \frac{\ln \frac{3}{2}}{2}
\;\leq\; \frac{1}{4}.
\ea
$$
Thus,
\beq\label{sum-Delta}
\ba{rcl}
\sum\limits_{i=0}^{k-1} \Delta_i &\refLE{sum-Delta-prel}& 
\frac{n L}{3 \mu} + \ln\Det(J_0^{-1}, J_k).
\ea
\eeq

Summing up \eqref{p-theta-prel} and using the fact that
$\psi_k \geq 0$, we obtain
\beq\label{p-theta}
\ba{rcl}
\frac{1}{7} \sum\limits_{i=0}^{k-1} p_i \theta_i^2 &\refLE
{p-theta-prel}& \psi_0 - \psi_k + \sum\limits_{i=0}^{k-1}
\Delta_i \;\leq\; \psi_0 + \sum\limits_{i=0}^{k-1} \Delta_i
\\
&\refEQ{sch-qn}& \psi(J_0, L B) + \sum\limits_{i=0}^{k-1}
\Delta_i \\
&\refEQ{def-psi}& \la J_0^{-1}, L B - J_0 \ra -
\ln\Det (J_0^{-1}, L B) + \sum\limits_{i=0}^{k-1} \Delta_i
\\
&\refLE{sum-Delta}& \la J_0^{-1}, L B - J_0 \ra + \frac{n L}
{3 \mu} - \ln \Det(J_k^{-1}, L B) \\
&\refLE{mu-L}& \la J_0^{-1}, \frac{L}{\mu} J_0 - J_0 \ra +
\frac{n L}{3 \mu} \\
&\refEQ{tr-det-A-inv-A}& n \left(\frac{L}{\mu} - 1 \right) +
\frac{n L}{3 \mu} \;\leq\; \frac{4}{3} \frac{n L}{\mu}.
\ea
\eeq
Since $(1 + t)^p \leq 1 + p t$ for all $t \geq -1$
and $0 \leq p \leq 1$, we further have
$$
\ba{rcl}
1 + \frac{M r_i}{2} &\leq& e^{\frac{M r_i}{2}}
\;\refLE{lam-lin}\; e^{\frac{M \lambda_0}{2}} \;\refLE
{lam-ini}\; \left( \frac{3}{2} \right)^{1/8} \\
&=& \sqrt{\left(\frac{3}{2} \right)^{1/4}}
\;\leq\; \sqrt{1 + \frac{1}{4} \cdot \frac{1}{2}} \;=\; 
\sqrt{\frac{9}{8}}.
\ea
$$
Therefore, by Lemma~\ref{lm-lam-prog} and the
arithmetic-geometric mean inequality,
$$
\ba{rcl}
\lambda_f(x_k) &\leq& \lambda_f(x_0) \prod\limits_{i=0}^
{k-1} \left[ \sqrt{\frac{9}{8}} \theta_i \right] \;=\; 
\frac{1}{\prod_{i=0}^{k-1} p_i^{1/2}} \left[ \left(
\frac{9}{8} \right)^k \prod\limits_{i=0}^{k-1} p_i
\theta_i^2 \right]^{1/2} \lambda_f(x_0) \\
&\leq& \frac{1}{\prod_{i=0}^{k-1} p_i} \left( \frac{9}{8}
\cdot \frac{1}{k} \sum\limits_{i=0}^{k-1} p_i \theta_i^2
\right)^{k/2} \lambda_f(x_0) \;\refLE{p-theta}\; \left(
\frac{9}{8} \cdot 7 \cdot \frac{4}{3} \frac{n L}{\mu k}
\right)^{k/2} \lambda_f(x_0) \\
&\leq& \left( \frac{21 n L}{2 \mu k} \right)^{k/2}
\lambda_f(x_0) \;\leq\; \left( \frac{11 n L}{\mu k} \right)^
{k/2} \lambda_f(x_0).$\qed$
\ea
$$

\section{Discussion}
\label{sec-disc}
\SetEQ

Let us compare the rates of superlinear convergence, that we
have obtained for the classical quasi-Newton methods, with
those of the greedy quasi-Newton methods
\cite{RodomanovNesterov2020}. For brevity, we discuss only
the BFGS method. Moreover, since the complexity bounds
for the general nonlinear case differ from those for the
quadratic one only in some absolute constants (both for the
classical and the greedy methods), we only consider the
case, when the objective function $f$ is quadratic.

As before, let $n$ be the dimension of the problem, $\mu$ be
the strong convexity parameter, $L$ be the Lipschitz
constant of the gradient of $f$, and $\lambda_f(x)$ be
the local norm of the gradient of $f$ at the point $x \in
\E$ (as defined by \eqref{quad-def-lam}). Also, let us
introduce the following condition number to simplify our
notation:
\beq\label{def-Q}
\ba{rcl}
Q &\Def& \frac{n L}{\mu} \;\geq\; 1.
\ea
\eeq

The greedy BFGS method \cite{{RodomanovNesterov2020}} is
essentially the classical BFGS algorithm (scheme
\eqref{quad-sch-qn} with $\phi_k \equiv 0$) with the only
difference that, at each iteration, the update direction
$u_k$ is chosen greedily according to the following rule:
$$
\ba{rcl}
u_k &\Def& \argmax\limits_{u \in \{e_1, \ldots, e_n\}} \frac{\la
G_k u, u \ra}{\la A u, u \ra},
\ea
$$
where $e_1, \ldots, e_n$ is a basis
in $\E$, such that $B^{-1} = \sum_{i=1}^n e_i e_i^*$.
For this method,
we have the following recurrence
(see \cite[Theorem~3.2]{RodomanovNesterov2020}):
$$
\ba{rcl}
\lambda_f(x_{k+1}) &\leq& \left( 1 - \frac{1}{Q} \right)^k Q
\lambda_f(x_k) \;\leq\; e^{-\frac{k}{Q}} Q \lambda_f(x_k),
\quad k \geq 0.
\ea
$$
Hence, its rate of superlinear convergence is described by
the expression
\beq\label{gr-rate}
\ba{rcl}
\lambda_f(x_k) &\leq& \lambda_f(x_0)
\prod\limits_{i=0}^{k-1} \left[ e^{-\frac{i}{Q}} Q \right]
\;=\; e^{-\frac{k (k-1)}{2 Q}} Q^k \lambda_f(x_0) \;\Def\;
A_k, \quad k \geq 0.
\ea
\eeq
Although the inequality \eqref{gr-rate} is valid for all $k
\geq 0$, it is useful only when
\beq\label{gr-sup-st}
\ba{rclrcl}
e^{-\frac{k (k-1)}{2 Q}} Q^k &\leq& 1 \qquad
\Longleftrightarrow \qquad k &\geq& 1 + 2 Q \ln Q.
\ea
\eeq
In other words, the relation \eqref{gr-sup-st} specifies the
moment, starting from which it becomes meaningful to speak
about the superlinear convergence of the greedy BFGS method.

For the classical BFGS method, we have the following bound 
(see~\eqref{quad-bfgs-sup}):
$$
\ba{rcl}
\lambda_f(x_k) &\leq& \left( \frac{Q}{k} \right)^{k/2}
\lambda_f(x_0) \;\Def\; B_k, \quad k \geq 1,
\ea
$$
and the starting moment of its superlinear convergence is
described as follows:
\beq\label{std-sup-st}
\ba{rclrcl}
\left( \frac{Q}{k} \right)^{k/2} &\leq& 1
\qquad \Longleftrightarrow \qquad
k &\geq& Q.
\ea
\eeq

Comparing \eqref{gr-sup-st} and \eqref{std-sup-st}, we see
that, for the standard BFGS, the superlinear convergence may
start slightly earlier than for the greedy one. However,
the difference is only in the logarithmic factor.

Nevertheless, let us show that, very soon after the
superlinear convergence of the greedy BFGS begins, namely,
after
\beq\label{def-K}
\ba{rcl}
K &\Def& 1 + 6 Q \ln(4 Q) \qquad (\;\refGE{def-Q}\; 7)
\ea
\eeq
iterations, it will be significantly faster than the
standard BFGS. Indeed,
\beq\label{AB-prel}
\ba{rcl}
\frac{A_k}{B_k} &=& e^{-\frac{k (k-1)}{2 Q}} Q^k \left( 
\frac{k}{Q} \right)^{k/2} \;=\; e^{-\frac{k (k-1)}{2 Q}} (Q
k)^{k/2} \\
&=& e^{-\frac{k (k-1)}{2 Q} + \frac{k}{2} \ln(Q k) } \;=\;
e^{-\frac{k (k-1)}{2 Q} \left[1 - \frac{Q \ln(Q k)} {k-1}
\right]}
\ea
\eeq
for all $k \geq 1$. Note that the function $t \mapsto
\frac{\ln t}{t}$ is decreasing on $[e, +\infty)$ (since its
logarithm $\ln \ln t - \ln t$ is a decreasing function of $u
= \ln t$ for $u \in [1, +\infty)$, which is easily verified
by differentiation). Hence, for all $k \geq K$, we have
(using first that $k \leq 2 (k-1)$ since $k \geq 2$)
$$
\ba{rcl}
\frac{Q \ln(Q k)}{k-1} &\leq& \frac{Q \ln(2 Q (k-1))}{k-1}
\;\leq\; \frac{Q \ln(2 Q (K-1))}{K-1} \;\refEQ{def-K}\; 
\frac{\ln\left( 12 Q^2 \ln(4 Q) \right)}{6 \ln(4 Q)} \\
&\leq& \frac{\ln(48 Q^3)}{6 \ln(4 Q)} \;\leq\; \frac{\ln(64
Q^3)}{6 \ln(4 Q)} \;=\; \frac{3 \ln(4 Q)}{6 \ln(4 Q)} \;=\; 
\frac{1}{2}.
\ea
$$
Consequently, for all $k \geq K$, we obtain
$$
\ba{rcl}
\frac{A_k}{B_k} &\refLE{AB-prel}& e^{-\frac{k (k-1)}{4 Q}}
\;\leq\; 1.
\ea
$$
Thus, after $K$ iterations, the rate of superlinear
convergence of the greedy BFGS is always better than that of
the standard BFGS. Moreover, as $k \to \infty$, the gap
between these two rates grows as $e^{-k^2/Q}$. At
the same time, the complexity of the Hessian update for the
greedy BFGS method is more expensive than for the standard
one.

\begin{acknowledgements}
The authors are thankful to two anonymous reviewers for
their valuable time and useful feedback.
\end{acknowledgements}

\section{Appendix}\label{sec-app}
\SetEQ

\BL\label{lm-proj}
Let $A, B : \E \to \E^*$ be self-adjoint linear operators
such that
\beq\label{proj-cond}
\ba{rclrcl}
0 &\prec& A &\preceq& B.
\ea
\eeq
Then, for any $u \in \E \setminus \{0\}$, we have
$$
\ba{rcl}
A - \frac{A u u^* A}{\la A u, u \ra} &\preceq& B - \frac{B u
u^* B}{\la B u, u \ra}.
\ea
$$
\EL

\proof
Indeed, for all $h \in \E$, we have
$$
\ba{rcl}
\la A h, h \ra - \frac{\la A u, h \ra^2}{\la A u, u \ra}
&=& \min\limits_{\alpha \in \R} \left[ \la A h, h \ra - 2
\alpha \la A h, u \ra + \alpha^2 \la A u, u \ra \right] \\
&=& \min\limits_{\alpha \in \R} \, \la A (h - \alpha u), h -
\alpha u \ra \\
&\refLE{proj-cond}& \min\limits_{\alpha \in \R} \, \la B (h
- \alpha u), h - \alpha u \ra \\
&=& \min\limits_{\alpha \in \R} \left[ \la B h, h \ra - 2
\alpha \la B h, u \ra + \alpha^2 \la B u, u \ra \right] \\
&=& \la B h, h \ra - \frac{\la B u, h \ra^2}{\la B u, u
\ra}.$\qed$
\ea
$$

\BL\label{lm-broyd-det-inv}
For any self-adjoint positive definite linear operators $A,
G : \E \to \E^*$, any scalar $\phi \in \R$, and any
direction $u \in \E \setminus \{0\}$, we have
\beq\label{det-Gp}
\ba{rcl}
\Det(G^{-1}, \Broyd_{\phi}(A, G, u)) &=& \phi \frac{\la A
G^{-1} A u, u \ra}{\la A u, u \ra} + (1-\phi) \frac{\la
A u, u \ra}{\la G u, u \ra}.
\ea
\eeq
\EL

\BR
Note that formula \eqref{det-Gp} is known in the literature
(see e.g. \cite[eq. (1.9)]{ByrdLiuNocedal1992}), although we
do not know any reference, which contains an explicit proof
of this result.
\ER

\proof
Denote $G_+ \Def \Broyd_{\phi}(A, G, u)$,
\beq\label{def-G0-s}
\ba{rclrcl}
G_0 &\Def& G - \frac{G u u^* G}{\la G u, u \ra} + \frac{A u
u^* A}{\la A u, u \ra},
\qquad
s &\Def& \frac{A u}{\la A u, u \ra} - \frac{G u}{\la G u, u
\ra}.
\ea
\eeq
Note that
$$
\ba{rcl}
G_+ &\refEQ{def-broyd}& G_0 + \phi \left[ \frac{\la G u, u
\ra A u u^* A}{\la A u, u \ra^2} + \frac{G u u^* G} {\la G
u, u \ra} - \frac{A u u^* G + G u u^* A}{\la A u, u \ra}
\right] \\
&=& G_0 + \phi \la G u, u \ra s s^*,
\ea
$$
and
\beq\label{su-0}
\ba{rcl}
\la s, u \ra &=& 0.
\ea
\eeq
Let $Q \Def G + \frac{A u u^* A}{\la A u, u \ra}$. Note that
\beq\label{Q-aux}
\ba{rclrcl}
Q u &=& G u + A u,
\qquad
Q G^{-1} A u &=& \left( 1 + \frac{\la A
G^{-1} A u, u \ra}{\la A u, u \ra} \right) A u,
\ea
\eeq
and $G_0 = Q - \frac{G u u^* G}{\la G u, u \ra}$. Therefore,
applying twice Lemma~\ref{lm-det-rank1}, we find that
$$
\ba{rcl}
\Det(G^{-1}, Q) &=& 1 + \frac{\la A G^{-1} A u, u \ra}{\la
A u, u \ra}, \\
\Det(Q^{-1}, G_0) &=& 1 - \frac{\la G Q^{-1}
G u, u \ra}{\la G u, u \ra} \;\refEQ{Q-aux}\; 1 - \frac{\la
G u - G Q^{-1} A u, u \ra}{\la G u, u \ra} \\
&=& \frac{\la G Q^{-1} A u, u \ra}{\la G u, u \ra} \;\refEQ
{Q-aux}\; \frac{\la A u, u \ra}{\la G u, u \ra \left ( 1 + 
\frac{\la A G^{-1} A u, u \ra}{\la A u, u \ra} \right)}.
\ea
$$
Hence,
\beq\label{det-G-inv-G0}
\ba{rcl}
\Det(G^{-1}, G_0) &\refEQ{det-mult}& \Det(G^{-1}, Q) \Det(Q^
{-1}, G_0) \;=\; \frac{\la A u, u \ra}{\la G u, u \ra}.
\ea
\eeq
Further, note that
\beq\label{G0-u-aux}
\ba{rclrcl}
G_0 u &\refEQ{def-G0-s}& A u,
\qquad
G_0 G^{-1} A u &\refEQ{def-G0-s}& \frac{\la A G^
{-1} A u, u \ra}{\la A u, u \ra} A u + A u - \frac{\la A
u, u \ra}{\la G u, u \ra} G u.
\ea
\eeq
So, applying Lemma~\ref{lm-det-rank1} again, we
obtain
\beq\label{det-G0-inv-Gp}
\ba{rcl}
\Det(G_0^{-1}, G_+) &\refEQ{def-G0-s}& 1 + \phi \la G u, u
\ra \la s, G_0^{-1} s \ra \\
&\refEQ{def-G0-s}& 1 + \phi \frac{\la G u, u \ra}{\la A u, u
\ra} \la s, G_0^{-1} A u - \frac{\la A u, u \ra}{\la G u, u
\ra} G_0^{-1} G u \ra \\
&\refEQ{G0-u-aux}& 1 + \phi \frac{\la G u, u \ra}{\la A u, u
\ra} \la s, G^{-1} A u - \frac{\la A G^{-1} A u, u \ra}{\la
A u, u \ra} G_0^{-1} A u \ra \\
&\refEQ{G0-u-aux}& 1 + \phi \frac{\la G u, u \ra}{\la
A u, u \ra} \la s, G^{-1} A u - \frac{\la A G^{-1} A u, u
\ra} {\la A u, u \ra} u \ra \\
&\refEQ{su-0}& 1 + \phi \frac{\la G u, u \ra}{\la A u, u
\ra} \la s, G^{-1} A u \ra \\
&\refEQ{def-G0-s}& 1 + \phi \frac{\la G u, u \ra}{\la A
u, u \ra} \la \frac{A u}{\la A u, u \ra} - \frac{G u}{\la G
u, u \ra}, G^{-1} A u \ra \\
&=& \phi \frac{\la G u, u \ra \la A G^{-1} A u, u \ra}{\la A
u, u \ra^2} + 1 - \phi.
\ea
\eeq
Consequently,
$$
\ba{rcl}
\Det(G^{-1}, G_+) &\refEQ{det-mult}& \Det(G^{-1}, G_0) \Det
(G_0^{-1}, G_+) \;\refEQ{det-G-inv-G0}\; \frac{\la A u, u
\ra}{\la G u, u \ra} \Det(G_0^{-1}, G_+) \\
&\refEQ{det-G0-inv-Gp}& \frac{\la A u, u \ra}{\la G u, u
\ra} \left( \phi \frac{\la G u, u \ra \la A G^
{-1} A u, u \ra}{\la A u, u \ra^2} + 1 - \phi \right) \\
&=& \phi \frac{\la A G^{-1} A u, u \ra}{\la A u, u \ra}
+ (1 - \phi) \frac{\la A u, u \ra}{\la G u, u \ra}.$\qed$
\ea
$$

\BL[Determinant of rank-1 perturbation]
\label{lm-det-rank1}
Let $A : \E \to \E^*$ be a self-adjoint positive definite
linear operator, $s \in \E^*$, and $\alpha \in \R$.
Then,
$$
\ba{rcl}
\Det(A^{-1}, A + \alpha s s^*) &=& 1 + \alpha \la s, A^{-1}
s \ra.
\ea
$$
\EL

\proof
Indeed, with respect to $A$, the operator $A + \alpha s s^*$
has $n-1$ unit eigenvalues and one eigenvalue $1 + \alpha
\la s, A^ {-1} s \ra$ (corresponding to the eigenvector
$A^{-1} s$).
\qed

\end{document}